\documentclass{amsart}

\usepackage{amssymb}

\newtheorem{thm}{Theorem}[section]
\newtheorem{prop}[thm]{Proposition}
\newtheorem{lem}[thm]{Lemma}
\newtheorem{cor}[thm]{Corollary}

\theoremstyle{definition}
\newtheorem{defi}[thm]{Definition}

\theoremstyle{remark}
\newtheorem{rem}[thm]{Remark}

\theoremstyle{example}
\newtheorem{ex}[thm]{Example}

\numberwithin{equation}{section}

\newcommand{\R}{\mathbb{R}^n}  
\newcommand{\Din}{D_{in}}
\newcommand{\Dout}{D_{out}}
\newcommand{\unv}{u_{\nu}}
\newcommand{\PR}{\mathcal{P}(\R)}
\newcommand{\nvGo}{\nu_{\Gamma_0}}
\newcommand{\SLam}{^\Lambda}
\newcommand{\EDt}{E^{\Delta t}}

\begin{document}

\title{Motion of sets by curvature and derivative of capacity potential}

\author{Hui Yu}
\address{Department of Mathematics, the University of Texas at Austin}
\email{hyu@math.utexas.edu}

\begin{abstract}
We study a geometric flow where the motion of a set is driven by the mean curvature of its boundary and the normal derivative of its capacity potential. We establish local well-posedness and propose two possible weak formulations that exist after singularities.
\end{abstract}

 \maketitle

\tableofcontents

\section{Introduction}
In this paper we study the following geometric evolution problem. Given two smooth open sets in $\R$, $D_{in}\subset\subset D_{out}$, and a smooth positive function $\phi$ on $ \overline{D_{in}}$, we study the evolution of a family of sets $\{\Omega_{t}\}_{t>0}$ that has initial configuration $D_{in}\subset\subset\Omega_0\subset\subset D_{out}$ and moves with outward normal velocity at $x\in\partial(\Omega_t)$ given by $$V(x,t)=u_{\nu}^2(x,t)-H(x,t).$$ Here $u(\cdot,t)$ is the capacity potential of $(\phi, \overline{\Din})$ relative to $\Omega_t$, namely, the minimizer of the Dirichlet energy over functions that coincide with $\phi$ on $\overline{D_{in}}$ and vanish outside $\Omega_t$. $u_{\nu}$ is its normal derivative, and $H(\cdot, t)$ is the mean curvature along $\Gamma_{t}=\partial \Omega_t$, taken to be positive for convex sets. Note here that we make no connectivity condition on $\Din$ and no boundedness condition on $\Dout$.

This can be seen as a parabolic version of a {\it elliptic problem} studied by Athanasoupoulos-Caffarelli-Kenig-Salsa \cite{ACKS}, where the authors investigated the existence and regularity of minimizers of the following energy $$\int|\nabla u|^2dx+Per(\{u>0\}),$$ over functions coinciding with $\phi$ on $\Din$ and vanishing outside $\Dout$. Here $Per$ is the perimeter of a set, to be understood in the sense of Giusti \cite{Giusti}. For its definition and properties the reader may refer to the wonderful exposition of Maggi \cite{Maggi}.  

Being a minimizer of the Dirichlet energy,  $u$ satisfies $\Delta u=0$ in $\{u>0\}$. The balancing condition on the free boundary $\partial\{u>0\}$, under our one-phase assumption, is $$\unv^2-H=0.$$   As a result, the free boundary in our problem, $\Gamma_t$, is driven towards the equilibrium of the area-Dirichlet integral. A related elliptic problem with more general types of bulk energies is studied by Mazzone \cite{Mazzone}. One should also consult Kim \cite{Kim} for more evolution laws  of this type.

Comparing with classic geometric evolutions, say, the mean curvature flow, one major difference lies in the dependence on the  the capacity potential, which in turn depends on the global shape of the domain. Our flow is thus a {\it nonlocal} curvature flow \cite{Card}\cite{CMP}. Such bulk interaction is absent in the mean curvature flow, which is completely described by a heat equation {\it on the interface}. Moreover, the geometry of our problem dictates that there would be no translation invariance nor interior-exterior reflection as in the mean curvature flow, which poses more difficulties.

If we again draw an analogy to the mean curvature flow, we see the study of geometric evolutions typically consists of two steps. Firstly, one establishes a well-posedness theory, starting from a nice geometric object as in Gage-Hamilton \cite{GH} or Evans-Spruck \cite{ES2}. After a certain lapse of time, however, singularities occur, and one faces the issue of defining a global continuation of the flow that is in some sense weaker, but still inherits certain key features of the original flow. 

For the mean curvature flow, three of such features are 1) the heat equation on the surface, 2) the inclusion principle and 3) the curve-shortening property. 

The heat equation leads to  the level-set formulation, first proposed by Osher-Sethian \cite{OS} and then analytically validated and extended by Chen-Giga-Goto \cite{CGG} and Evans-Spruck \cite{ES1}. One might argue it also leads to the varifold formulation by Brakke \cite{Brakke}, although in a less direct fashion. 
The inclusion principle gives rise to weak formulations that come from geometric comparisons as in the works of De Giorgi \cite{DeGiorgi} and Barles-Souganidis \cite{BS}. And the curve-shortening property inspired the formulation through minimizing movements by Almgren-Taylor-Wang \cite{ATW} and Luckhaus-Sturzenhecker \cite{LS}. If one accepts that the perimeter functional is the limit of Ginzburg-Landau \cite{Modica}, then it is also responsible for the formulation via singular limits of diffusion equations as by Bronsard-Kohn \cite{BK} and Evans-Soner-Souganidis \cite{ESS}. Of course the literature on the mean curvature flow and its weak formulations is particularly rich. The reader should refer to lecture notes of Souganidis \cite{Souganidis} for a more thorough reference list and the book by Bellettini \cite{Bellettini} for many modern developments. 

We follow these two steps for our problem. 

The first step is to show the local existence and uniqueness of a smooth flow, starting from a very nice geometry. Given a smooth initial configuration $\Omega_0$ and $\Gamma_0=\partial\Omega_0$, the problem is essentially to study the following coupled systems of PDEs
$$\begin{cases}\Delta_x u(\cdot, t)=0 & \text{ in $\Omega_t\backslash\overline{\Din}$,}\\
u(\cdot, t)=\phi &\text{ on $\overline{\Din}$,} \\ u(\cdot, t)=0 &\text{ on $\Gamma_t=\partial\Omega_t$},
\end{cases}$$ and $$\begin{cases}V(\cdot,t)=\unv^2(\cdot, t)-H_{\Gamma_t}(\cdot, t) &\text{ on $\Gamma_t$,}\\ \Gamma_0=\Omega_0.\end{cases}$$ Following traditions, the first system is referred to as the {\it bulk} equation, the second as the equation of motion. 

Here we follow the original idea of Hanzawa \cite{Hanzawa} by first translating the problem into the fixed domain $\Omega_0$ via a change of variable, and then reducing the problem to a parabolic operator on the deformation, which is  solved using a fixed point argument very much inspired by the work of Chen-Reitich \cite{CR}. It is noteworthy that one might also apply the theory of maximal regularity \cite{Angenent} as in Pr\"uss-Saal-Simonett \cite{PSS}. Since we are only concerned with smooth initial configurations, these two theories give similar results, namely:

\begin{thm}
Given a smooth open set $\Omega_0$ with $\Din\subset\subset\Omega_0\subset\subset\Dout$, there exists some $T>0$ and a unique smooth flow starting from $\Omega_0$ on time interval $[0,T]$.
\end{thm} 

The definition of a smooth flow to our problem is given in the next section. We'd like to mention that the positivity of $\phi$ does not play a role in this local well-posedness theory, and that we could as well consider a two-phase problem. Also, when $\phi=0$, our surface moves by the classical mean curvature flow, and $\Din$ effectively becomes an obstacle \cite{ACN}. 

Since we are considering very general domains $\Din$, $\Dout$ and boundary data $\phi$, such a well-posedness theory is only available over very short time intervals.  Consider the very simple case where $\Din$ consists of two balls, $B_1$ and $B_2$, with unit radii, centered at two points of distance $1000$ to each other. Let $\phi=1$ and $\Dout$ be a very large ball.  Note that the optimal configuration for the elliptic problem consists of two slightly larger balls $\tilde{B}_1\supset B_1$ and $\tilde{B}_2\supset B_2$. As a result, if we take as our initial configuration a simply connected domain by connecting $\tilde{B}_1$ to $\tilde{B}_2$ with  a very thin {\it neck}, then the flow drives the domain into two separate pieces, each containing one of the two optimal balls.  It is readily seen that here the topology of the domain undergoes a transition from a simply connected set to two disjoint sets, resulting in a breaking of any possible smooth solutions. Note that in this work we do not provide estimate on the lifespan of smooth flows as in Evans-Spruck \cite{ES3} or Giga-Yama-Uchi \cite{GY}. Needless to say, such estimates would be very interesting and useful. However, they seem difficult to get for the very general situation considered here. 

After the short time of smooth flows, we propose two possible formulations for a weak flow that exists globally in time. Due to the nonlocal nature of our problem, it is not likely that one can reduce it to the study of just a surface equation or the motion of level sets of some function. It is also difficult to see how it can be reached as the limit of diffusion equations. For that one would need some scaling relation that keeps both the bulk equation and the curvature in the limit.  This leaves us with the options of a variational approach and a geometric-comparison approach. 

The idea of a variational approach to curvature flows originated from the work of Almgren-Taylor-Wang \cite{ATW} and  Luckhaus-Sturzenhecker \cite{LS}. Take again the mean curvature flow as an example, they study the following type of minimization problem $$E\mapsto Per(E)+\frac{1}{\Delta t}\int_{E\Delta E_0}dist(x,\partial E_0)dx.$$ Here $\Delta t>0$ is the discrete time step, $E_0$ the initial configuration, $E\Delta E_0$ the symmetric difference between the two sets. The motivation is that the mean curvature flow is the gradient flow of  the perimeter functional, and the energy above is a discretized version of the gradient flow.  Moreover, the Euler-Lagrange equation  for this minimization problem is $$H_{\partial E}(x)=\frac{dist(x,\partial E_0)}{\Delta t} \text{ for $x\in\partial E$ }.$$ Taking the right-hand side as the discretized velocity at which $x\in\partial E$ is moving away from $\partial E_0$, one sees the relation to the mean curvature flow. 

In Almgren-Taylor-Wang, they studied the discrete motion generated by this minimizing movement scheme, $\{E(k\Delta t)\}_{k\in\mathbb{N}}$, where $E(0)=E_0$, and $E(k\Delta t)$ is a minimizer of the same minimization problem with $E_0$ replaced by $E((k-1)\Delta t)$. Using various tools from geometric measure theory, they show uniform compactness of this scheme as $\Delta t\to 0$, the limiting motion from which is termed {\it a flat flow}. More technical is that under various conditions this flat flow agrees with the smooth flow as long as the latter exists. It can also be shown that along this flat flow the perimeter decreases \cite{Chambolle}, which further justifies the flat flow as a weak continuation of the mean curvature flow.

For our problem, it is natural to study the following minimization procedure $$(u, E)\mapsto \int|\nabla u|^2dx+ Per(E)+\frac{1}{\Delta t}\int_{E\Delta E_0}dist(x,\partial E_0)dx,$$ where the tuple $(u, E)$ is taken over certain class of admissible functions and sets to be clarified later. Here we are relying on the heuristic that our flow drives the sets to optimal shapes for the area-Dirichlet integral considered by Athanaspoulos-Caffarelli-Kenig-Salsa \cite{ACKS}.  Following Algrem-Taylor-Wang we are aiming at certain compactness of the discrete motion as $\Delta t\to 0$. The nonlocal Dirichlet energy, however, poses challenging problems. Fortunately, there is a very rich theory on shape optimization problems involving the Dirichlet energy \cite{AC}\cite{BB}. Combining this theory with the estimates in Almgren-Taylor-Wang, we successfully establish a H\"older continuous flat flow that exists globally in time:

\begin{thm}
Starting from any set of finite perimeter $E_0$ with $\Din\subset\subset E_0\subset\subset\Dout$ with $\Dout$ bounded, there exists a flat flow for our problem. Moreover, this flow is H\"older continuous in time. 
\end{thm} 

The definition of a flat flow is given later, as well as the meaning of H\"older continuity in this context. 

We do not have uniqueness of flat flows, which is somehow expected since nor is there a uniqueness theory for the elliptic problem, and our flat flow in a sense is built from the elliptic problem. However, this does pose serious challenges for the consistency of our flat flow with the smooth flow, which seems very difficult for the general situation that we are considering here. This lack of uniqueness also leads to the lack of a semigroup property for flat flows. 

These problems can be easily tackled, however, if one uses our second weak formulation based on geometric comparison principles. As for the case of the mean curvature flow, this is based on the inclusion principle, namely, if we have two initial domains $U\subset V$, then the classical flows starting from these two initial data remain ordered.  It is the idea of De Giorgi \cite{DeGiorgi} to use smooth flows as test flows, and to call any flow that expands faster than all smooth flows a {\it barrier}. One sees here the analogy with the theory of viscosity solutions to elliptic/ parabolic PDEs. The least supersolution in that theory translates to the minimal barrier, which is another natural weak continuation of the mean curvature flow.

For a nonlocal curvature flow, one generally does not have such nice property as the inclusion principle, which is in a sense a {\it local} property. Luckily for us, our flow does enjoy this property. Consider two smooth domains $U\subset V$ at time zero. Suppose, on the contrary, that $V_t$ failed to contain $U_t$ at some later time, then there would be some {\it critial} time, at which $U_t$ is still in $V_t$ but $\partial U_t$ would touch $\partial V_t$ at some point $x$. However, since our domains are still ordered at this time, their respective capacity potentials satisfy $0\le u\le v$. With $u(x)=v(x)=0$ one has $\unv^2\le v_{\nu}^2$. Also, the inclusion of sets also implies $H_{\partial U_t}(x)\ge H_{\partial V_t}(x)$, hence the outward velocity of $U$ is less than the velocity of $V$ at $x$, and $U$ would stay contained in $V$. This heuristic is justified as in \begin{thm}
Suppose $f,g:[a,b]\to\PR$ are two smooth flows.

If $f(a)\subset g(a)$ then $f(b)\subset g(b)$.
\end{thm}

The local well-posedness theory, this observation of inclusion principle and the idea of De Giorgi naturally lead to the second weak formulation using minimal barriers. This formulation enjoys many nice properties. Of particular interest to us is the following:
\begin{thm}
The flow of minimal barriers enjoys the semigroup property, and is consistent with the smooth flow as long as the latter exists.
\end{thm} 
 
 We also give a conditional result concerning the long-term behaviour  of the flow of minimal barriers:
 
 \begin{thm}
If there are enough smooth flows converging to the optimal configuration of the elliptic problem, then so does the flow of minimal barriers.
\end{thm}  

See the last section for the precise statement.
 
 Geometric comparison properties have also proven to be very useful in the study of higher regularity as in Caffarelli-Salsa \cite{CS}. However, we do not pursuit this direction in this work.

This paper is organized as follows: In the next section we give the formal definition of a smooth flow, and establish a local well-posedness theory starting from nice configurations. After this we turn to the two weak formulations. In the third section we give the definition of a flat flow following the original work of Almgren-Taylor-Wang \cite{ATW}. Then we give several geometric-measure-theoretic estimates before we prove the compactness of the discrete flow, which gives the existence of a flat flow for our problem. In the fourth section we deal with the formulation by minimal barriers. This is done by first prove an inclusion principle for smooth flows, which is necessary for the comparison. Then we establish certain geometric properties of the flow of minimal barriers, in particular, the uniqueness, the semigroup property and the consistency with the smooth flow. 

As already mentioned, this work leaves many interesting problems open. For instance, some estimate on the lifespan of smooth flows would be both interesting and useful. One could try to show the consistency of the three flows, at least for short time. This would in particular imply a semigroup property for flat flows as well as an energy decreasing property for the flow of minimal barriers. One might also study regularity properties of the two weak flows, especially the minimal barrier flow. Another interesting question is the long-term behaviour of our flows. Do they converge to optimal configurations of the elliptic problem, if, say, we start from a configuration close to one?  These problems seem difficult in the very general situation considered here, but one could first gain some intuition by computations with special geometries, for instance, in lower dimensions or using radial configurations.

\section{Local well-posedness}
In this section and the rest, $\R$ denotes the $n$-dimensional Euclidean space. $\PR$ is the collection of subsets of $\R$. 

For $x\in\R$ and $E\in\PR$, the distance function is defined by $$dist_{E}(x)=dist(x,E)=\inf_{y\in E}|x-y|.$$For $E,F\in\PR$, the distance function between them is defined by $$dist(E,F)=\inf_{x\in E, y\in F}|x-y|.$$ The signed distance is defined by $$d_E(x)=dist_E(x)-dist_{\R\backslash E}(x).$$ Under this convention, the signed distance is nonpositive inside the set and nonnegative outside.

For a smooth $\Omega\in\PR$ and $x\in\partial\Omega$, $\nu_{\partial\Omega}(x)$ is the outward unit normal at $x$.  Note that $\nu_{\partial\Omega}(x)=\nabla d_\Omega(x)$. The mean curvature at this point of $\partial\Omega$ is then $$H_{\partial\Omega}(x)=\Delta d_{\Omega}(x).$$ Note that in particular the mean curvature of a convex set is nonnegative. For a function $u$ that is smooth on $\Omega$, $\unv(x)=u_{\nu_{\partial\Omega}}(x)$ is the normal derivative {\it from inside} at $x$, i.e., $$\unv(x):=\lim_{t\to 0^{-}}\frac{u(x+t\nu(x))-u(x)}{t}.$$Since we are only concerned with the quantity $\unv^2$, the sign of this normal derivative is not a matter for us. However, it is noteworthy that we only need information {\it inside} the domain where $u$ is smooth in defining this normal derivative.

We now steal from Bellettini \cite{Bellettini} the definition of a smooth flow:

\begin{defi}$f$ is a smooth flow if 
\begin{enumerate}
\item there exist $a<b$ in $\mathbb{R}$ such that $f:[a,b]\to\PR$,
\item for each $t\in [a,b]$, $f(t)$ is closed, and
\item if we denote by $d(\cdot,t)=d_{f(t)}(\cdot)$ for $t\in [a,b]$, then for each $t$ there is an open set $A_t\supset\partial f(t)$ such that $d(\cdot,\cdot)$ is smooth in $\bigcup_{[a,b]}A_t\times\{t\}.$
\end{enumerate}
\end{defi} 
\begin{rem}
The choice of closed sets $f(t)$ over open sets has no significance on the theory to be developed. 
\end{rem}

We give the definition of {\it a smooth flow driven by mean curvature and the normal derivative of capacity potential}, which is to be called {\it a smooth MCND flow}. This definition is again modelled after the definition of the smooth mean curvature flow in \cite{Bellettini}.

\begin{defi}
A smooth MCND flow is a tuple $(u,f;[a,b])$ where $u:\R\times [a,b]\to\mathbb{R}$ and $f$ is a smooth flow on $[a,b]$. We further require that 
\begin{enumerate}
\item $\Din\subset\subset f(t)\subset\subset\Dout$ for all $t\in[a,b]$,
\item $u(x,t)=\phi(x)$ for all $x\in\overline{\Din}$ and $u(x,t)=0$ if $x\notin f(t)$ for all $t\in[a,b]$,
\item $u$ satisfies at each time $t\in[a,b]$ $$\Delta_x u(\cdot,t)=0 \text{ in $Int(f(t))\backslash\overline{\Din}$},$$
\item for $t\in[a,b]$ and $x\in\partial f(t)$, the signed distance function to $f(t)$ satisfies
$$\frac{\partial d}{\partial t}(x,t)=\Delta_xd(x,t)-\unv^2(x,t).$$
\end{enumerate}
\end{defi} 

\begin{rem}Here and afterwards, a subscript $x$ indicates the differential operator only acts on spatial variables. 

The first three items are to guarantee that $u(\cdot, t)$ is the capacity potential of $(\phi,\overline{\Din})$ relative to $f(t)$. The last item ensures that the outward normal velocity of the boundary is $\unv^2(x,t)-H(x,t)$.
\end{rem} 

\begin{rem}
In this section the exterior domain $\Dout$ plays no role at all. It may as well be $\R$.
\end{rem}

\begin{rem}
Note that the potential $u$ is completely determined by the flow $f$ via the Dirichlet problem, hence we would sometimes just use $f$ to denote a smooth MCND flow. Also note that when $\partial f(t)$ is smooth, $u$ is smooth in the spatial variables inside $f(t)$ as a consequence of elliptic regularization \cite{GT}. Therefore all quantities in the definition are well-defined.
\end{rem} 

\subsection{Notations and preparations}

Now let $\Omega_0$ be a smooth open set with $\Din\subset\subset\Omega_0\subset\subset\Dout$. $\Gamma_0=\partial \Omega_0$. The goal is to show the existence and uniqueness of a smooth MCND flow $f:[0,T]\to\PR$ with $Int(f(0))=\Omega_0$ for a possibly small $T>0$. As already mentioned in the introduction, what follows is very much modelled on the argument as in Chen-Reitich \cite{CR}.

Note that under our assumption, $\Gamma_0$ is a smooth $(n-1)$-dimensional manifold in $\R$ without boundary (not necessarily connected).  

Let $s'=(s_1,s_2,\dots,s_{n-1})\in U\subset\mathbb{R}^{n-1}$ denote a generic coordinate of a point in $\Gamma_0$, and let $X^0:U\to X^0(U)\subset\Gamma_0$ be a generic coordinate system of $\Gamma_0$.

For $L_0>0$ depending only on $\Din, \Dout$ and $\Gamma_0$, the map $X:U\times [-L_0,L_0]\to \R$ defined by $$(s_1,s_2,\dots,s_{n-1},s_n)\mapsto X^0(s_1,s_2,\dots,s_{n-1})+s_n\nu_{\Gamma_0}(X^0(s_1,s_2,\dots,s_{n-1}))$$ is smooth. Moreover, it is a diffeomorphism between its domain and its range, $N$, which is a neighborhood of $\Gamma_0$. We can further assume $\Din \subset\subset N\subset\subset\Dout$. 

Let $S=(S',S^n)=(S^1,S^2,\dots,S^{n-1},S^n)$ denote the inverse of $X$. Note that $S'$ maps a point to the coordinates of its projection onto $\Gamma_0$. $S^n$ is just the signed distance from $\Gamma_0$.

We will also need the derivatives of $S$ with respect to $x$, which is a matrix, $\frac{\partial S}{\partial x}$, the $i$-th column of which is $$\alpha^i(S(x))=\begin{pmatrix}
\frac{\partial S^i}{\partial x_1}\\\frac{\partial S^i}{\partial x_2}\\ \cdot\\ \cdot\\ \frac{\partial S^i}{\partial x_n}
\end{pmatrix}(x).$$ The second order derivatives of $S^i$ with respect to $x$ are denoted by matrix $A^{i}$ whose $(s,k)$-entry is $\frac{\partial^2 S^i}{\partial x_s\partial x_k}$. Note that these matrices only depend on the initial configuration $\Gamma_0$.

Now let $\Lambda:U\times [0,T]\to\mathbb{R}$ be some function with $\Lambda(s',0)=0$ for all $s'\in U$ and $|\Lambda|\le L_0$. This function denotes the {\it deformation} at time $t$ from our initial configuration. In particular, for each fixed $t$, if we define $$\Gamma_t:=\{x+\Lambda(S'(x),t)\nu_{\Gamma_0}(x):x\in\Gamma_0\},$$ then $\Lambda$ induces a natural diffeomorphism between $\Gamma_0$ to $\Gamma_t$ $$\theta_{\Lambda}(x,t):=x+\Lambda(S'(x),t)\nu_{\Gamma_0}(x).$$ 

It would be convenient to have a diffeomorphism of the entire Euclidean space that agrees with $\theta_{\Lambda}(\cdot, t)$ on $\Gamma_0$. To this end, we pick $\zeta:\mathbb{R}\to\mathbb{R}$, a smooth nonnegative function which vanishes outside $[-L_0,L_0]$ and stays $1$ inside $[-\frac{3}{4}L_0,\frac{3}{4}L_0]$. Allow me to use the same $\theta_{\Lambda}(\cdot, t)$ to denote the following diffeomorphism from $\R$ to $\R$:$$x\mapsto x+\zeta(S^n(x))\Lambda(S'(x),t)\nvGo(X^0(S'(x))) \text{ if $x\in N$,}$$ and $$x\mapsto x \text{ otherwise}.$$

Define $\Omega_t=\theta_{\Gamma}(\Omega_0,t)$, then it's clear that $\Gamma_t=\partial\Omega_t$. Also $\Din\subset\subset\Omega_t\subset\subset\Dout$.

We now define for $x\in N$ $$\Phi_{\Lambda}(x,t)=S^n(x)-\Lambda(S'(x),t),$$ which can be extended to the entire $\R$ easily. This is the so-called {\it defining function} of $\Gamma_t$ in the sense that $$\Gamma_t=\{\Phi_{\Lambda}(\cdot,t)=0\}.$$

As a consequence, the outward normal velocity at $x\in\Gamma_t$ is 
\begin{align*}
V(x,t)&=\frac{-\frac{\partial}{\partial t}\Phi(x,t)}{|\nabla_x\Phi(x,t)|}\\&=\frac{\frac{\partial}{\partial t}\Lambda(S'(x),t)}{|\nabla_x\Phi(x,t)|}.
\end{align*}
The mean curvature is 
\begin{align*}
H(x,t)=&div_x(\frac{\nabla_x\Phi}{|\nabla_x\Phi|})(x,t)\\=&-\frac{1}{|\nabla_x\Phi|}(\Sigma_{i,j=1}^{n-1}a_{ij}(S'(x),\Lambda(S'(x),t),\nabla_{s'}\Lambda(S'(x),t))\frac{\partial^2\Lambda}{\partial s_i\partial s_j}\\&  +b(S'(x),\Lambda(S'(x),t),\nabla_{s'}\Lambda(S'(x),t)).
\end{align*}

Here $$a_{ij}(s',s_n,p^1,\dots,p^{n-1})=\alpha^i\cdot\alpha^j-\frac{\Sigma_{k,l=1}^{n-1}(p^k\alpha^k\cdot\alpha^i)(p^l\alpha^l\cdot\alpha^j)}{1+|\Sigma_{k=1}^{n-1}p^k\alpha^k|^2},$$ and 

\begin{align*}
b(s',s_n,p^1,\dots, p^{n-1})=&\Sigma_{k=1}^{n-1}p^ktrace(A^k)-trace(A^n)\\&-\frac{\Sigma_{k,l=1}^{n-1}p^kp^l(A^n\alpha^k)\cdot\alpha^l}{1+|\Sigma_{k=1}^{n-1}p^k\alpha^k|^2}+\frac{\Sigma_{k,l,m=1}^{n-1}p^kp^lp^m(A^m\alpha^k)\cdot\alpha^l}{1+|\Sigma_{k=1}^{n-1}p^k\alpha^k|^2}.
\end{align*}
For details of this computation the reader should consult Chen-Reitich \cite{CR}. To simplify notations we would write $$a_{ij}^\Lambda(s',t)=a_{ij}(s',\Lambda(s',t),\nabla_{s'}\Lambda(s',t))$$ and $$b^\Lambda(s',t)=b(s',\Lambda(s',t),\nabla_{s'}\Lambda(s',t)).$$

Consequently, if we find a function $u:\R \times [0,T]\to\mathbb{R}$ such that for each $t$,
\begin{equation}\begin{cases}\Delta_xu(x,t)=0 &\text{ in $\Omega_t\backslash\overline{\Din}$,}\\ u(x,t)=\phi &\text{ in $\overline{\Din}$,}\\ u(x,t)=0 &\text{ on $\Gamma_t$},\end{cases}\end{equation}
and a smooth function $\Lambda$ such that 
\begin{equation}\begin{cases}\frac{\partial}{\partial t}\Lambda(S'(x),t)=\Sigma_{ij=1}^{n-1}a^\Lambda_{ij}\frac{\partial^2\Lambda}{\partial s_i\partial s_j}(S'(x),t)+b^\Lambda+|\nabla_x\Phi(x,t)|\unv^2(x,t) &\text{ in $t\in [0,T]$, $x\in\Gamma_t$,}\\ \Lambda(s',0)=0 \text{ for all $s'$},\end{cases}\end{equation}
 then $f(t):=\overline{\Omega_t}$ gives a smooth MCND flow in the sense of Definition 2.3.

The last bit of preparation we need is some parabolic function spaces. Let $M$ be a generic subset in $\R$, which could be a domain in $\R$ or some hypersuface. For each $\alpha, \beta\ge 0$, the H\"older space $C^{\beta,\alpha}(M\times [0,T])$ is the completion of the space of smooth functions on that domain under the norm $$\|f\|_{\beta,\alpha}=\|f\|_{\mathcal{L}^{\infty}}+\sup_{m\in M}\|f(m,\cdot)\|_{\alpha}+\sup_{t\in[0,T]}\|f(\cdot,t)\|_{\beta},$$ where $\|g\|_{\gamma}$ is the usual $\gamma$-H\"older norm. 

The parabolic H\"older space $C^{\beta,\frac{\beta}{2}}(M\times [0,T])$ is the completion of the space of smooth functions under the norm $$\|f\|_{\beta,\frac{\beta}{2}}=\Sigma_{0\le|\alpha|+2k\le\beta} \|D^{\alpha}_mD_t^kf\|_{\beta,0}+\Sigma_{0\le|\alpha|+2k\le\beta} \|D^{\alpha}_mD_t^kf\|_{0,\beta/2},$$here $\alpha$ is a multi-index. 

The reader might refer to Friedman \cite{Friedman} or Lady\u zenskaya-Solonnikov-Ural'ceva \cite{LSU} for more properties of these spaces.  What is useful for us is that there is a constant $C=C(\Gamma_0)$ such that $$\|a_{ij}^\Lambda\|_{(k-1)+\alpha,\frac{(k-1)+\alpha}{2}}\le C\|\Lambda\|_{k+\alpha,\frac{k+\alpha}{2}}$$ and $$\|b^\Lambda\|_{(k-1)+\alpha,\frac{(k-1)+\alpha}{2}}\le C\|\Lambda\|_{k+\alpha,\frac{k+\alpha}{2}}$$ for $k=1,2$.

\subsection{Reduction to a fixed domain} We follow the idea of Hanzawa \cite{Hanzawa} to reduce our problem to a fixed domain. Define a new function $$v(x,t)=u(y,t),$$ where $y=\theta_{\Lambda}(x,t)$.

If $u$ is a solution to (2.1) then $v$ solves the following 
\begin{equation}
\begin{cases}
div_x(M^\Lambda(x,t)\nabla v(x,t))=0 &\text{ in $\Omega_0\backslash\overline{\Din}$,}\\v(x,t)=\phi &\text{ in $\overline{\Din}$,}\\ v(x,t)=0 &\text{ along $\partial \Omega_0$},
\end{cases}
\end{equation} where $$M^\Lambda(x,t)=det(D_x\theta_{\Lambda}(x,t))(D_x\theta_{\Lambda}(x,t))^{T}.$$ Here $D_x$ is the differential of a map. And $T$ denotes matrix transposition. 

We now turn to the equation of motion (2.2). 

Since $\Gamma_t$ is the zero level surface of $u(\cdot, t)$, one has for $y\in\Gamma_t$, 
\begin{align*}
|\unv(y,t)|=&|\nabla u(y,t)|\\=&|(D_y\theta_{\Lambda}^{-1}(y,t))^{T}\nabla_xv(x,t)|.
\end{align*}

Now a direct computation gives $$(D\theta_{\Lambda}^{-1})^T(y,t)=Id-\nu_{\Gamma_0}(S'(y))\otimes\nabla_y\Lambda(S'(y),t)-\Lambda(S'(y),t)(D_y\nu_{\Gamma_0}(S'(y)))^T,$$ where $\otimes$ denotes the tensor product between two matrices. 

Considering the fact that $\Gamma_0$ is the zero level surface of $v$, one has that $\nabla_xv$ is parallel to $\nu_{\Gamma_0}$. Also, $\nabla_y\Lambda(S'(y),t)$ is perpendicular to $\nu_{\Gamma_0}$, and $\nabla_xv$ is in the kernel of  $D_y\nu_{\Gamma_0}(S'(y))$. Consequently, $$|\nabla u(y,t)|\\=|(D_y\theta_{\Lambda}^{-1}(y,t))^{T}\nabla_xv(x,t)|=(1+|\nabla_y\Lambda(S'(y),t)|^2)^{1/2}|\nabla_xv(x,t)|.$$

If we let $$c^{\Lambda}=|\nabla_y\Phi_{\Lambda}(y,t)|(1+|\nabla_y\Lambda(S'(y),t)|^2)^{1/2},$$ the equation of motion (2.2) becomes a parabolic equation on the boundary-less manifold $\Gamma_0$:
\begin{equation}
\begin{cases}\frac{\partial}{\partial t}\Lambda(s',t)=\Sigma_{ij}a^\Lambda_{ij}(s',t)\frac{\partial^2\Lambda}{\partial s_i\partial s_j}(s',t)+b(s',t)+c^\Lambda(s',t)|\nabla_xv(X^0(s'),t)|^2 &\text{ in $\Gamma_0\times(0,T]$,}\\ \Lambda(s',0)=0&\text{ in $\Gamma_0$}.\end{cases}
\end{equation} 

Note that we have reduced the systems (2.1) and (2.2) to fixed domains, and it is clear that solving (2.3) and (2.4) is equivalent to solving the original systems after a change of variable. In the next subsection, we solve these systems by a fixed point argument. 

\subsection{A fixed point argument}The main difficulty now is the coupling between $\Lambda$ and $v$. All the coefficients of the equations depend on $\Lambda$ while the right-hand side of the equation for $\Lambda$ depends on $v$.  To solve this problem, we first fix the coefficients at a certain $\Lambda$, then solve for $v$. Plug this $v$ into the parabolic equation with coefficients corresponding to the same $\Lambda$, we obtain a solution $\lambda$. 

This procedure gives rise to a map $\Lambda\mapsto\lambda$, which can be shown to be a contraction mapping on certain subset of some parabolic H\"older space. Then the well-posedness follows directly from the contraction mapping principle of Banach \cite{Banach}. 

Throughout this subsection, $\alpha$ is some number in $(0,1).$

We first solve the bulk equation:\begin{lem}
Given $\Lambda\in C^{2+\alpha,\frac{2+\alpha}{2}}(\Gamma_0\times[0,T])$, there is a unique solution to

\begin{equation*}
\begin{cases}
div_x(M^\Lambda(x,t)\nabla v^\Lambda(x,t))=0 &\text{ in $\Omega_0\backslash\overline{\Din}$}\\v^\Lambda(x,t)=\phi &\text{ in $\overline{\Din}$,}\\ v^{\Lambda}(x,t)=0 &\text{ along $\partial \Omega_0$}.
\end{cases}\end{equation*} 
Furthermore, for some $C=C(\Din,\Omega_0,\|\Lambda\|_{2+\alpha,\frac{2+\alpha}{2}})$ we have the following estimates
$$\|v^\Lambda\|_{2+\alpha,0}\le C(1+\|\phi\|_{C^3})$$ and 
$$\|\nabla v^\Lambda\|_{\frac{1}{2}(1+\alpha),\frac{1}{4}(1+\alpha)}\le C(1+\|\phi\|_{C^3}).$$
\end{lem} 

\begin{rem}
The dependence on $C^3$-norm of $\phi$ is far from sharp. However, this suffices our purpose. 
\end{rem} 
\begin{proof}
From its definition one has $$\|M\SLam\|_{1+\alpha,\frac{1+\alpha}{2}}\le C(n,\|\Lambda\|_{2+\alpha,\frac{2+\alpha}{2}}).$$ Thus theory of elliptic equations of divergence type \cite{GT} gives the existence and uniqueness of the solution $v\SLam$. 

Moreover, $$\|v\SLam\|_{2+\alpha,0}\le C(1+\|\phi\|_{C^3})$$ for some $C=C(\Din,\Omega_0,\|\Lambda\|_{2+\alpha,\frac{2+\alpha}{2}})$.

We now trade some spatial regularity for temporal regularity. 

For $t_1\neq t_2$ in $[0,T]$,$$div_x(M\SLam(x,t_1)\nabla_xv\SLam(x,t_1))=0=div_x(M\SLam(x,t_2)\nabla_xv\SLam(x,t_2)) \text{ in $\Omega_0\backslash\overline{\Din}$}.$$

Therefore, $$\begin{cases}
div_x(M\SLam(t_1)\nabla_x(v\SLam(t_1)-v\SLam(t_2))=div_x((M\SLam(t_2)-M\SLam(t_1))\nabla_xv\SLam(t_2)) &\text{ in $\Omega_0\backslash\overline{\Din}$,}\\v\SLam(t_1)-v\SLam(t_2)=0 &\text{ on $\partial(\Omega_0\backslash\overline{\Din})$}.
\end{cases}$$

Apply elliptic regularization to $v(t_1)-v(t_2)$, one thus obtains

\begin{align*}\|v\SLam(t_1)-v\SLam(t_2)\|_{C^{1+\frac{1+\alpha}{2}}}\le &C\|(M\SLam(t_2)-M\SLam(t_1))\nabla_xv\SLam(t_2)\|_{C^{\frac{1+\alpha}{2}}}\\ \le &C\|(M\SLam(t_2)-M\SLam(t_1))\|_{C^{\frac{1+\alpha}{2}}}\|\nabla_xv\SLam(t_2)\|_{C^{\frac{1+\alpha}{2}}}\\\le&C\|M\SLam\|_{1+\alpha,\frac{1+\alpha}{2}}|t_1-t_2|^{\frac{1+\alpha-(1+\alpha)/2}{2}}\|v\SLam(t_2)\|_{C^{2+\alpha}}\\\le&C(1+\|\phi\|_{C^3})|t_1-t_2|^{\frac{1+\alpha}{4}}.\end{align*} Note that we have used the parabolic estimate on $M\SLam$ as well as the spatial regularity on $v\SLam(t_2)$. Obviously this leads to the second estimate in the lemma. 
\end{proof} 

We now turn to the equation of motion:
\begin{lem}
Given $\Lambda\in C^{2+\alpha,\frac{2+\alpha}{2}}$, let $v\SLam$ be as in the previous lemma. Then there is a unique solution $\lambda\SLam$ to the following

\begin{equation*}
\begin{cases}\frac{\partial}{\partial t}\lambda\SLam(s',t)=\Sigma_{ij}a\SLam_{ij}(s',t)\frac{\partial^2\lambda\SLam}{\partial s_i\partial s_j}(s',t)+b(s',t)+c^\Lambda(s',t)|\nabla_xv\SLam(X^0(s'),t)|^2 &\text{ in $\Gamma_0\times(0,T]$,}\\ \lambda\SLam(s',0)=0&\text{ in $\Gamma_0$}.\end{cases}
\end{equation*} 

Furthermore, for any $m_0>0$, there is  $T_0=T_0(\phi, m_0)>0$ such that $\|\Lambda\|_{2+\alpha,\frac{2+\alpha}{2}}\le m_0$ and $T<T_0$ implies $$\|\lambda\SLam\|_{C^{2+\alpha,\frac{2+\alpha}{2}}}\le m_0.$$
\end{lem} 

\begin{proof}

For $\Lambda\in C^{2+\alpha,\frac{2+\alpha}{2}}$, it is obvious $a^{\Lambda}_{ij},b\SLam$ and $c\SLam$ are in $C^{1+\alpha,\frac{1+\alpha}{2}}$. Meanwhile, the second estimate in the previous lemma implies $\nabla_xv\SLam\in C^{\frac{1+\alpha}{2},\frac{1+\alpha}{4}}$. Thus the standard theory of parabolic equations implies the well-posedness of the problem \cite{Friedman}\cite{LSU}. 

As for the estimate, we first note that $$\|a\SLam_{ij}\|_{\frac{1+\alpha}{2},\frac{1+\alpha}{4}},\|b\SLam\|_{\frac{1+\alpha}{2},\frac{1+\alpha}{4}},\|c\SLam\|_{\frac{1+\alpha}{2},\frac{1+\alpha}{4}}\le C(m_0)$$if $\|\Lambda\|_{2+\alpha,\frac{2+\alpha}{2}}\le m_0$.

Also the estimate from the previous lemma implies $$\|\nabla v\SLam\|_{\frac{1+\alpha}{2},\frac{1+\alpha}{4}}\le C(m_0)(1+\|\phi\|_{C^3}).$$

Combining these estimates with parabolic regularization, one has 
$$\|\lambda\SLam\|_{2+\frac{1+\alpha}{2},1+\frac{1+\alpha}{4}}\le C(m_0,\phi).$$

With the initial data $\lambda\SLam(\cdot, 0)=0$, we have the following interpolation 
\begin{align*}
\|\lambda\SLam\|_{C^{2+\alpha,\frac{2+\alpha}{2}}}&\le T^{\frac{1-\alpha}{4}}\|\lambda\SLam\|_{2+\frac{1+\alpha}{2},1+\frac{1+\alpha}{4}}\\&\le C(m_0,\phi)T^{\frac{1-\alpha}{4}}\\&\le m_0
\end{align*}
once we choose $T$ small. 
\end{proof} 

The previous lemma implies that the map $\Lambda\mapsto\lambda\SLam$ is an endomorphism on the space $\{f\in C^{2+\alpha,\frac{2+\alpha}{2}}:\|f\|_{C^{2+\alpha,\frac{2+\alpha}{2}}}\le m_0\}$ for any $m_0>0$ once we shorten the time interval enough. 

Our next goal is to show that the map is contracting. 

\begin{lem}
Given $\Lambda^1,\Lambda^2$ with $\|\Lambda^j\|_{2+\alpha,\frac{2+\alpha}{2}}\le m_0$ $j=1,2$, and let $\lambda^1$ and $\lambda^2$ be solutions to the equation of motion as in the previous lemma with $\Lambda$ replaced by $\Lambda^1$ and $\Lambda^2$ respectively. 

If we take $T$ small enough, depending on $\Din, \phi, \Omega_0$ and $m_0$,  then $$\|\lambda^1-\lambda^2\|_{2+\alpha,\frac{2+\alpha}{2}}\le\frac{1}{2}\|\Lambda^1-\Lambda^2\|_{2+\alpha,\frac{2+\alpha}{2}}.$$
\end{lem} 
\begin{proof}
Let $v_j$ be the solution to the Dirichlet problem in Lemma 2.7 when $\Lambda$ is replaced with $\Lambda^j$, $j=1,2$. Let $d=v_1-v_2$, then $d$ satisfies
$$\begin{cases}
div(M^{\Lambda^1}\nabla d)=div((M^{\Lambda^2}-M^{\Lambda^1})\nabla v_2) &\text{ in $\Omega_0\backslash\overline{\Din}$,}\\d=0 &\text{ on $\partial(\Omega_0\backslash\Din)$}.
\end{cases}$$

Note that by definition of $M\SLam$, $$\|M^{\Lambda^1}-M^{\Lambda^2}\|_{1+\alpha,\frac{1+\alpha}{2}}\le C\|\Lambda^1-\Lambda^2\|_{2+\alpha,\frac{2+\alpha}{2}},$$ where $C$ depends on $m_0, \Din,$ and $ \Omega_0$.

Also, invoking the first estimate in Lemma 2.7, $$\|v_2\|_{2+\alpha,0}\le C(1+\|\phi\|_{C^3}).$$

Combining these two estimate, the following is a consequence of elliptic regularization
\begin{align*}
\|d\|_{2+\alpha,0}\le &C(\|M^{\Lambda^1}\|_{1+\alpha,0})\|(M^{\Lambda^2}-M^{\Lambda^1})\nabla v_2\|_{1+\alpha,0}\\\le &C(m_0)\|M^{\Lambda^1}-M^{\Lambda^2}\|_{1+\alpha,0}\|\nabla v_2\|_{1+\alpha,0}\\\le&C(1+\|\phi\|_{C^3})\|\Lambda^1-\Lambda^2\|_{2+\alpha,\frac{2+\alpha}{2}}.
\end{align*}Here the constant $C$ depends on $m_0, \Din, \Omega_0$. 

We now turn to the temporal regularity of $d$. 

Again for $t_1\neq t_2$ in $[0,T]$,
\begin{align*}
div(M^{\Lambda^1(t_1)}\nabla(d(t_1)-d(t_2)))=&div((M^{\Lambda^1(t_2)}-M^{\Lambda^1(t_1)})\nabla d(t_2))\\+div(&(M^{\Lambda^2(t_2)}-M^{\Lambda^1(t_2)})(\nabla v_2(t_1)-\nabla v_2(t_2)))\\+div&([(M^{\Lambda^2(t_1)}-M^{\Lambda^1(t_1)})-(M^{\Lambda^2(t_2)}-M^{\Lambda^1(t_2)})]\nabla v_2(t_2)).
\end{align*}

We now estimate each of the terms on the right-hand side.

For the first term, we have $$\|\nabla d(t_2)\|_{1+\alpha,0}\le C(1+\|\phi\|_{C^3})\|\Lambda^1-\Lambda^2\|_{2+\alpha,\frac{2+\alpha}{2}},$$ and 
\begin{align*}
\|M^{\Lambda^1(t_1)}-M^{\Lambda^1(t_2)}\|_{\frac{1+\alpha}{2},0}\le &|t_1-t_2|^{\frac{1+\alpha}{4}}\|M^{\Lambda^1}\|_{1+\alpha,\frac{1+\alpha}{2}}\\ \le &C(m_0)|t_1-t_2|^{\frac{1+\alpha}{4}}.
\end{align*} Consequently, $$\|(M^{\Lambda^1(t_2)}-M^{\Lambda^1(t_1)})\nabla d(t_2)\|_{\frac{1+\alpha}{2},0}\le C(1+\|\phi\|_{C^3})\|\Lambda^1-\Lambda^2\|_{2+\alpha,\frac{2+\alpha}{2}}|t_1-t_2|^{\frac{1+\alpha}{4}}.$$

For the second term, we first invoke the second estimate in Lemma 2.7 to get 
$$\|\nabla v_2(t_1)-\nabla v_2(t_2)\|_{\frac{1+\alpha}{2},0}\le C(1+\|\phi\|_{C^3})|t_1-t_2|^{\frac{1+\alpha}{4}}.$$ Meanwhile, it's obvious 
$$\|M^{\Lambda^2(t_2)}-M^{\Lambda^1(t_2)}\|_{\frac{1+\alpha}{2},0}\le C\|\Lambda^2-\Lambda^1\|_{2+\alpha,\frac{2+\alpha}{2}}.$$
As a result, 
$$\|(M^{\Lambda^2(t_2)}-M^{\Lambda^1(t_2)})(\nabla v_2(t_1)-\nabla v_2(t_2))\|_{\frac{1+\alpha}{2},0}\le C(1+\|\phi\|_{C^3})\|\Lambda^2-\Lambda^1\|_{2+\alpha,\frac{2+\alpha}{2}}|t_1-t_2|^{\frac{1+\alpha}{4}}.$$

For the last term, first note, as in Lemma 2.7,  $$\|\nabla v_2(t_2)\|_{1+\alpha,0}\le C(1+\|\phi\|_{C^3}).$$

Also, expand the definition of $M\SLam$, note that determinant is a smooth function, and that all the matrices involved are of norm less than some number bounded by $m_0$, we have the following
$$\|(M^{\Lambda^2(t_1)}-M^{\Lambda^1(t_1)})-(M^{\Lambda^2(t_2)}-M^{\Lambda^1(t_2)})\|_{\frac{1+\alpha}{2},0}\le C\|\Lambda^2-\Lambda^1\|_{2+\alpha,\frac{2+\alpha}{2}}|t_1-t_2|^{\frac{1+\alpha}{4}}.$$
Therefore we have the following on the last term on the right-hand side 
\begin{align*}\|[(M^{\Lambda^2(t_1)}-M^{\Lambda^1(t_1)})-(M^{\Lambda^2(t_2)}-M^{\Lambda^1(t_2)})]&\nabla v_2(t_2)\|_{\frac{1+\alpha}{2},0}\\\le &C(1+\|\phi\|_{C^3})\|\Lambda^2-\Lambda^1\|_{2+\alpha,\frac{2+\alpha}{2}}|t_1-t_2|^{\frac{1+\alpha}{4}}.\end{align*}

These estimates give the following via elliptic regularization
$$\|\nabla d(t_1)-\nabla d(t_2)\|_{\frac{1+\alpha}{2},0}\le C(1+\|\phi\|_{C^3})\|\Lambda^2-\Lambda^1\|_{2+\alpha,\frac{2+\alpha}{2}}|t_1-t_2|^{\frac{1+\alpha}{4}}.$$That is, 
$$\|\nabla d\|_{\frac{1+\alpha}{2},\frac{1+\alpha}{4}}\le C(1+\|\phi\|_{C^3})\|\Lambda^2-\Lambda^1\|_{2+\alpha,\frac{2+\alpha}{2}}.$$

Now the equation for $\lambda^1-\lambda^2$ is 
$$\frac{\partial}{\partial t}(\lambda^1-\lambda^2)-a_{ij}^{\Lambda^1}\frac{\partial^2(\lambda^1-\lambda^2)}{\partial s_i\partial s_j}=(a^{\Lambda^1}_{ij}-a^{\Lambda^2}_{ij})\frac{\partial^2\lambda^2}{\partial s_i\partial s_j}+b^{\Lambda^1}-b^{\Lambda^2}+c^{\Lambda^1}|\nabla v_1|^2-c^{\Lambda^2}|\nabla v_2|^2.$$

Again we can estimate each of the terms on the right-hand side, by invoking the estimate on $\lambda^2$ as in Lemma 2.9, and the temporal estimate on $\nabla d$ that we just obtained. The method is similar to the previous step so we omit the details. In the end, one has that the $C^{\frac{1+\alpha}{2},\frac{1+\alpha}{4}}$-norm of the right-hand side is bounded by a constant multiple of $\|\Lambda^1-\Lambda^2\|_{2+\alpha,\frac{2+\alpha}{2}}.$ Then parabolic regularization gives 
$$\|\lambda^1-\lambda^2\|_{2+\frac{1+\alpha}{2},1+\frac{1+\alpha}{4}}\le C\|\Lambda^1-\Lambda^2\|_{2+\alpha,\frac{2+\alpha}{2}}.$$

Therefore the same interpolation argument as in the previous lemma gives $$\|\lambda^1-\lambda^2\|_{2+\alpha,\frac{2+\alpha}{2}}\le\frac{1}{2}\|\Lambda^1-\Lambda^2\|_{2+\alpha,\frac{2+\alpha}{2}}$$ if one chooses $T$ small enough.
\end{proof} 

Now we can finally give the proof of Theorem 1.1.

\begin{proof}(of Theorem 1.1.)
The map $\Lambda\mapsto\lambda$ is a contraction on $\{f\in C^{2+\alpha,\frac{2+\alpha}{2}}:\|f\|_{C^{2+\alpha,\frac{2+\alpha}{2}}}\le m_0\}$ for any $m_0>0$ once $T$ is small enough. Banach's contraction principle then gives a unique fixed point of this map. It is clear such a fixed point $\Lambda$ with the corresponding $v^{\Lambda}$ would solve (2.3) and (2.4). After a change of variable the $\Lambda$ and $u$ would solve (2.1) and (2.2). Moreover, a simple bootstrapping would imply smoothness of $\Lambda$. As commented after Equation (2.2), $f(t)=\overline{\Omega_t}$ is a smooth MCND flow with $Int(f(0))=\Omega_0$.
\end{proof} 

\begin{rem}
It seems interesting that the positivity of $\phi$ does not play a role in the well-posedness theory. Hence we actually have such a theory for the two-phase problem. Also it is interesting to note that when $\phi=0$, our surface moves by its mean curvature. In this sense, our problem is a generalization of the mean curvature flow. The set $\Din$ then becomes effectively an {\it  obstacle} for the mean curvature flow. See Almeida-Chambolle-Novaga \cite{ACN} for a related problem.
\end{rem}

\subsection{A radial example}
As mentioned in the introduction, it would be very interesting to obtain estimates on the life span of smooth flows. However, this seems rather challenging for the general situation considered here. In particular, we posed no topological restriction on $\Din, \Omega_0$ at all. 

The following very easy example is to illustrate that under certain geometries one expects a global smooth flow, that converges to the optimal configuration for the area-Dirichlet integral.

\begin{ex}
For this simple example we take $\Din=B_1(0)$, $\Dout=\R$ and $\phi=1$. The choice of $\Dout$ is for convenience only, some very large ball would also suffice. 

By a simple comparison argument we see the unique minimizer to the area-Dirichlet integral is a radial function of the form  $$u_{R}(x)=\begin{cases}
\frac{R^{n-2}}{R^{n-2}-1}(\frac{1}{|x|^{n-2}}-\frac{1}{R^{n-2}}) &\text {in $B_R\backslash B_1$,}\\1&\text{in $B_1$}\\ 0&\text{ outside $B_R$}.
\end{cases}$$

Then to optimize the area-Dirichlet integral, we just need to minimize 
\begin{align*}
\int_{B_R\backslash B_1}|\nabla u_R|^2dx+Per(B_R)&=\int_{\partial B_1}\unv dH^{n-1}+n\omega_nR^{n-1}\\&=(n-2)\frac{R^{n-2}}{R^{n-2}-1}n\omega_n+n\omega_nR^{n-1},
\end{align*}here $\omega_n$ is the volume of the unit ball in $\R$. 

Then it reduces to a problem in calculus, and we see that for the minimizer $R_{opt}$ solves $$R_{opt}(R_{opt}^{n-2}-1)^2=(n-2)^2/(n-1).$$

We now turn to our geometric flow. 

For simplicity, we take, as the initial configuration, $f(0)=\overline{B_{L_0}(0)}$. Then the radial symmetry reduces our problem to the following evolution on the length of radius if we let $f(t)=\overline{B_{L(t)}(0)}$:
$$\begin{cases}
L(0)=L_0 &\text{ $t=0$,}\\ \frac{d}{dt}L(t)= (\frac{n-2}{L(L^{n-2}-1)})^2-\frac{n-1}{L} &\text{ $t>0$.}\end{cases}$$

Note now that $R_{opt}>1$ and $\frac{d}{dt}L>0$ if $L<R_{opt}$, we see that $L$ is driven away from the sigularity $L=1$.  In particular, this flow exists globally in time. 

Also, the right-hand side of the equation for $L$ is uniformly strictly decreasing with respect to $L$ when $L$ is closed to $R_{opt}$, and it stays uniformly away from $0$ when $L$ is away from $R_{opt}$, we see $L\to R_{opt}$ as $t\to\infty$ . 

\end{ex}

\section{The variational approach}
The idea of the variational approach comes from Almgren-Taylor-Wang \cite{ATW} and Luckhaus-Sturzenhecker \cite{LS}, which gives a weak continuation of the smooth mean curvature flow that inherits its gradient flow structure. Since the stationary point of our flow is also a minimizer for some optimization problem, it is natural to propose a weak formulation for this MCND flow via a minimizing movement scheme \cite{DeGiorgi2}.

\subsection{Approximate flows and flat flows} 
In this subsection we give the definition of flat flows, which was invented by Almgren-Taylor-Wang to study geometric flows from a variational perspective. 

But before that, let us first recall that the data for our problem consist of two smooth open sets $\Din$ and $\Dout$ in $\R$, with the former compactly contained in the latter, and a smooth positive function $\phi$ on $\overline{\Din}$. 

To avoid certain triviality, we assume in this section $n\ge 2$. 

Without loss of generality, let's assume $0\in\Din$. Being smooth, we might assume $\Din$ satisfies the interior ball condition with constant $\theta>0$. That is, for any point $x\in\partial\Din$, we can find a ball $B_r(z)$ inside $\Din$ such that $x\in\partial B_r(z)$ and $r\ge\theta$. 

In this section and this section alone we assume that  $\Dout$ is bounded. Denote $D:=\sup_{z\in \Dout}|z|$.

Fix some positive time step $\Delta t>0$, we now give the definition of {\it a $\Delta t$-approximate flow}.

Given a set of finite perimeter \cite{Giusti}\cite{Maggi} $E_0$ as our initial configuration, we look at the following minimization problem \begin{equation}
(u,E)\mapsto \mathcal{E}(u,E):=\int|\nabla u|^2dx+Per(E)+\frac{1}{\Delta t}\int_{E\Delta E_0}dist(x,\partial E_0)dx,
\end{equation} where $(u, E)$ is taken from the following admissible family 
$$\mathcal{A}:=\{(u,E):u\in H^1_0(\Dout), u=\phi \text{ on $\overline{\Din}$}, \Din\subset E\subset\Dout, \{u> 0\}\subset E \text{ a.e.}\}.$$

\begin{rem}
In Almgren-Taylor-Wang \cite{ATW} the minimization problem is considered in the entire $\R$, which requires a bound on the diameter of possible minimizers. This is done in Almgren-Taylor-Wang using the fact that by cutting a set with hyperplanes we can reduce its perimeter. This shows any possible minimizer is contained in the convex hull of the initial configuration, and in particular, the diameter of any possible minimizer is less than the diameter of the initial configuration. 

This is no longer the case for our problem. Therefore we {\it artificially} assume that everything lies inside a bounded set $\Dout$, which has the effect of making the last term in (3.1) a lower order term. 

We also remark that in two dimensions, we can get a natural bound on the diameter of minimizers and get rid of the assumption on $\Dout$. Such a bound can be obtained in the following way: the energy is less than some number $C=C(E_0)$, and in particular $Per(E)\le C$. Then we note that the perimeter of the convex hull of $E$ has less perimeter than $E$, and that the diameter of a convex set in $\mathbb{R}^2$ is bounded by its perimeter. Thus we have a bound of the diameter of $E$ in terms of $E_0$. Thus all possible minimizers are contained inside a ball, and there is no need to impose the exterior domain $\Dout.$
\end{rem} 

We return to the definition of a $\Delta t$-approximate flow.  

Once we find a minimizer $E^{\Delta t}(\Delta t)$ of (3.1), we can repeat the same procedure with $E^{\Delta t}(\Delta t)$ taking the place of $E_0$. In this way, one obtains a sequence of sets of finite perimeter $\{E^{\Delta t}(k\Delta t)\}_{k\in\mathbb{N}}$, which is a discretized version of our flow. 

\begin{defi}
A $\Delta t$-approximate flow starting from $E_0$ is a map $E^{\Delta t}:[0,\infty)\to\PR$ such that 
\begin{enumerate}
\item $E^{\Delta t}(0)=E_0$,
\item $E^{\Delta t}(k\Delta t)$ is a minimizer of the energy with initial configuration $E^{\Delta t}((k-1)\Delta t)$, and \item $\EDt(s\Delta t)=\EDt(k\Delta t)$ for $k\le s<(k+1)$. 
\end{enumerate}
\end{defi}  

\begin{rem}
Since the last term in (3.1) is of lower order, the existence of these minimizers follows relatively easily from the elliptic theory \cite{ACKS}. Also, as in the elliptic theory, we have a drastic lack of uniqueness. 
\end{rem} 

A {\it flat flow} is a limit of these discrete flows as $\Delta t\to 0$:

\begin{defi}
A map $E^0:[0,\infty)\to\PR$ is a flat flow starting from $E_0$ if $E^0(0)=E_0$ and there is a subsequence of $\Delta t$-approximate flows starting from $E_0$ with $$|E(t)\Delta\EDt(t)|\to 0$$ locally uniformly in time as $\Delta t\to 0$.
\end{defi}

The H\"older continuity of a flat flow is defined in the same $\mathcal{L}^1$ topology:
\begin{defi}
Such a flow is $\alpha$-H\"older continuous if there is some constant $C$ such that $$|E^{0}(s)\Delta E^0(t)|\le C|t-s|^{\alpha}.$$
\end{defi} 

\subsection{The minimization}
In this subsection we collect several results about the minimization problem (3.1) that are useful in constructing a flat flow. Most of the arguments are variants of ones in Athanasopoulos-Caffarelli-Kenig-Salsa \cite{ACKS} and Mazzone \cite{Mazzone}.

\begin{prop}
Given $E_0$, there exists a minimizer $(u,E)$.
\end{prop} 
\begin{proof}
The energy is obviously nonnegative. Taking a minimizing sequence $(u_n,E_n)$, then there are uniform bounds on the perimeters of $E_n$ and the Dirichlet energy of $u_n$. 

Compactness of sets of finite perimeter and Sobolev embedding then give a pair $(u,E)$ such that $$u_n\to u \text{ weakly in $H^1_0$}, $$ and $$E_n\to E \text{ in $\mathcal{L}^1$}.$$ The semicontinuity of the Dirichlet energy and the perimeter gives $$\int|\nabla u|^2dx\le\liminf \int|\nabla u_n|^2dx,$$and $$Per(E)\le\liminf Per(E_n).$$

Also Fatou's gives $$\frac{1}{\Delta t}\int_{E\Delta E_0}dist(x,\partial E_0)dx\le\liminf\frac{1}{\Delta t}\int_{E_n\Delta E_0}dist(x,\partial E_0)dx.$$

Consequently $\mathcal{E}(u,E)\le\liminf\mathcal{E}(u_n,E_n)$.

Moreover, it is also clear that all the criteria for admissibility are respected by weak $H^1_0$-convergence and $\mathcal{L}^1$-convergence. We see $(u,E)\in\mathcal{A}$, and it is a minimizer. 
\end{proof} 

\begin{prop}
For a minimizer, $0\le u\le\sup\phi$.
\end{prop} 
\begin{proof}
Let $\overline{u}:=\min\{u,\sup\phi\}$. Then $(\overline{u},E)\in\mathcal{A}$.

As a result, 
\begin{align*}
0\le&\int|\nabla\overline{u}|^2dx-\int|\nabla u|^2dx\\ =&-\int_{\{u>\sup\phi\}}|\nabla u|^2dx.
\end{align*}

Thus $|\{u>\sup\phi\}|=0$.

Similarly, if we take $\underline{u}=\max\{u,0\}$, then $(\underline{u},E)\in\mathcal{A}$.

\begin{align*}
0\le&\int|\nabla\underline{u}|^2dx-\int|\nabla u|^2dx\\ =&-\int_{\{u<0\}}|\nabla u|^2dx.
\end{align*}

Thus 
$|\{u<0\}|=0.$
\end{proof} 

\begin{prop}
Any minimizer satisfies $\Delta u\ge 0$ in $\Dout\backslash\overline{\Din}$.
\end{prop} 
\begin{proof}
Let $\eta\le 0$ be a smooth function compactly supported outside $\overline{\Din}$. Then for $t>0$ the pair $(u+t\eta,E)$ is admissible. 

Optimality of $u$ gives 
\begin{align*}
0\le&\int|\nabla(u+t\eta)|^2dx-\int|\nabla u|^2dx\\ =&2t\int\nabla u\cdot\nabla\eta dx+t^2\int|\nabla\eta|^2dx.
\end{align*}Dividing both sides by $t$ and sending $t\to0$ gives $$0\le\int\nabla u\cdot\nabla\eta dx.$$
\end{proof} 

Our next proposition says that a minimizer $E$ stays at a distance to $\Din$.

\begin{prop}
Let $E$ be a minimizer, then $$dist(\partial E,\partial\Din)\ge \min\{C(n,\theta,D)\sqrt{\Delta t}, (D-\theta), dist(\partial\Dout, \partial\Din)\}.$$
\end{prop} 
\begin{proof}
Suppose $$d_0=dist(\partial E,\partial\Din)=|x_0-y_0|$$where $x_0\in\partial E$ and $y_0\in\Din$. 

We assume $d_0<(D-\theta)$ and $d_0<dist(\partial\Dout,\partial\Din)$ and prove $$d_0\ge C(n,\theta,D)\sqrt{\Delta t}.$$

Let $B_\theta(z_0)\subset\Din$ be a ball tangent to $\partial\Din$ at $y_0$.

For $\epsilon>0$ small, define a function $v$ by 
$$\begin{cases}
\Delta v=0 &\text{ in $B_{\theta+d_0+\epsilon}(z_0)\backslash B_\theta(z_0)$,}\\v=\min\phi &\text{ in $\overline{B_\theta(z_0)}$,}\\ v=0 &\text{ along $\partial B_{\theta+d_0+\epsilon}(z_0)$.}
\end{cases}$$
Since $d_0<(D-\theta)$ and $d_0<dist(\partial\Dout,\partial\Din)$, all these perturbations stay in $\Dout$ if $\epsilon$ is small. 

Build a competitor by defining $\tilde{u}=\max\{u,v\}$ and $\tilde{E}=E\cup B_{\theta+d_0+\epsilon}(z_0)$.

Note that $\Din\subset\tilde{E}\subset\Dout$ by definition. Also $\tilde{u}\in H^1_0(\Dout)$ agrees with $\phi$ on $\Din$ since there $v\le\min\phi\le u$. In addition, $\{\tilde{u}>0\}\subset\tilde{E}$. As a result, $(\tilde{u},\tilde{E})\in\mathcal{A}.$ 

We compare each of the terms in (3.1).
\begin{equation*}
Per(\tilde{E})-Per(E)=H^{n-1}(\partial B_{\theta+d_0+\epsilon}(z_0)\cap E^c)-Per(E;B_{\theta+d_0+\epsilon}(z_0)).
\end{equation*}

To estimate the right-hand side, we first note that, if we let $\lambda=\frac{1}{\theta+d_0+\epsilon}$ be the nonzero eigenvalue of $D^2d_{B_{\theta+d_0+\epsilon}(z_0)}$ along $\partial B_{\theta+d_0+\epsilon}(z_0)$, then  inside $B_{\theta+d_0+\epsilon}(z_0)$ the nonzero eigenvalue is $\frac{\lambda}{1+d_{B_{\theta+d_0+\epsilon}(z_0)}\lambda}.$

Note that in $B_{\theta+d_0+\epsilon}(z_0)\backslash E$, $d_{B_{\theta+d_0+\epsilon}(z_0)}\ge -d_0-\epsilon$ since $B_{\theta}(z_0)\subset E$. Thus $$\frac{\lambda}{1+d_{B_{\theta+d_0+\epsilon}(z_0)}\lambda}\le\frac{\lambda}{1-(d_0+\epsilon)\lambda}=\frac{1}{\theta}.$$ As a result, $$\Delta d_{B_{\theta+d_0+\epsilon}(z_0)}\le \frac{n-1}{\theta}$$ in $B_{\theta+d_0+\epsilon}(z_0)\backslash E.$

Consequently, $$\int_{B_{\theta+d_0+\epsilon}(z_0)\backslash E}\Delta d_{B_{\theta+d_0+\epsilon}(z_0)}dx\le\frac{n-1}{\theta}|B_{\theta+d_0+\epsilon}(z_0)\backslash E|.$$ 

On the other hand \begin{align*}\int_{B_{\theta+d_0+\epsilon}(z_0)\backslash E}\Delta d_{B_{\theta+d_0+\epsilon}(z_0)}dx&=\int_{\partial B_{\theta+d_0+\epsilon}(z_0)\backslash E}\nabla d_{B_{\theta+d_0+\epsilon}(z_0)}\cdot\nu dH^{n-1}\\+&\int_{\partial E\cap B_{\theta+d_0+\epsilon}(z_0)}\nabla d_{B_{\theta+d_0+\epsilon}(z_0)}\cdot\nu dH^{n-1}\\&\ge \int_{\partial B_{\theta+d_0+\epsilon}(z_0)\backslash E} dH^{n-1}-\int_{\partial E\cap B_{\theta+d_0+\epsilon}(z_0)}dH^{n-1}\\&=H^{n-1}(\partial B_{\theta+d_0+\epsilon}(z_0)\cap E^c)-Per(E;B_{\theta+d_0+\epsilon}(z_0)).\end{align*} Here the boundary of a set as well as the outward normal derivatives are to be understood in the sense of Giusti \cite{Giusti} or Maggi \cite{Maggi}. Also note that we used crucially the fact that $\nabla d_{B_{\theta+d_0+\epsilon}(z_0)}=\nu$ on $\partial B_{\theta+d_0+\epsilon}(z_0)$, and that $|\nabla d_{B_{\theta+d_0+\epsilon}(z_0)}\cdot\nu|\le1$. 

Combining these estimates, one has $$Per(\tilde{E})-Per(E)\le \frac{n-1}{\theta}|B_{\theta+d_0+\epsilon}(z_0)\backslash E|.$$

Meanwhile, the distance integral can be easily estimated
$$\frac{1}{\Delta t}\int_{\tilde{E}\Delta E_0}dist(x,\partial E_0)dx-\frac{1}{\Delta t}\int_{E\Delta E_0}dist(x,\partial E_0)dx\le\frac{D}{\Delta t}|B_{\theta+d_0+\epsilon}(z_0)\backslash E|.$$ Note that here the boundedness of $\Dout$ plays a definitive role. 

Finally we turn to the Dirichlet energy.
Since $u=\tilde{u}$ inside $\Din$ and outside $B_{\theta+d_0+\epsilon}(z_0)$, one has $$\int|\nabla u|^2dx-\int|\nabla \tilde{u}|^2dx=\int_{B_{\theta+d_0+\epsilon}(z_0)\backslash\Din}|\nabla u|^2-|\nabla\tilde{u}|^2dx.$$
Since on $\partial(B_{\theta+d_0+\epsilon}(z_0)\backslash\Din)$ $u=\tilde{u}$, inside $B_{\theta+d_0+\epsilon}(z_0)\backslash\Din$ $u-\tilde{u}\le 0$ and $\Delta\tilde{u}\ge 0$, one has $$\int_{B_{\theta+d_0+\epsilon}(z_0)\backslash\Din}\nabla\tilde{u}\cdot\nabla(u-\tilde{u})\ge 0.$$ Thus
\begin{align*}
\int_{B_{\theta+d_0+\epsilon}(z_0)\backslash\Din}|\nabla u|^2-|\nabla\tilde{u}|^2dx&\ge\int_{B_{\theta+d_0+\epsilon}(z_0)\backslash\Din}|\nabla u|^2-|\nabla\tilde{u}|^2-2\nabla\tilde{u}\cdot\nabla(u-\tilde{u})dx\\&\ge\int_{B_{\theta+d_0+\epsilon}(z_0)\backslash E}|\nabla u|^2-|\nabla\tilde{u}|^2-2\nabla\tilde{u}\cdot\nabla(u-\tilde{u})dx\\&=\int_{B_{\theta+d_0+\epsilon}(z_0)\backslash E}|\nabla v|^2dx.
\end{align*}
Here the second inequality follows from the convexity of the map $p\mapsto |p|^2$, which leads to $|\nabla u|^2-|\nabla\tilde{u}|^2-2\nabla\tilde{u}\cdot\nabla(u-\tilde{u})\ge0$. The last line follows from the definition of $\tilde{u}$.

Since $v$ is the radial harmonic function in $B_{\theta+d_0+\epsilon}(z_0)\backslash B_{\theta}(z_0)$, we invoke the direct estimate $$|\nabla v|\ge\min\phi((\theta+d_0+\epsilon)^{n-1}(\frac{1}{\theta^{n-2}}-\frac{1}{(\theta+d_0+\epsilon)^{n-2}}))^{-1}.$$

Therefore \begin{align*}\int|\nabla u|^2dx-\int|\nabla \tilde{u}|^2dx&\ge\\(\min\phi)^2((\theta+d_0+\epsilon)^{n-1}&(\frac{1}{\theta^{n-2}}-\frac{1}{(\theta+d_0+\epsilon)^{n-2}}))^{-2}|B_{\theta+d_0+\epsilon}(z_0)\backslash E|.\end{align*}

The optimality of $(u,E)$ then implies \begin{align*}\frac{n-1}{\theta}|B_{\theta+d_0+\epsilon}(z_0)\backslash E|+\frac{D}{\Delta t}|B_{\theta+d_0+\epsilon}(z_0)\backslash E|&\ge\\(\min\phi)^2((\theta+d_0+\epsilon)^{n-1}(\frac{1}{\theta^{n-2}}-&\frac{1}{(\theta+d_0+\epsilon)^{n-2}}))^{-2}|B_{\theta+d_0+\epsilon}(z_0)\backslash E|.\end{align*}
Divide by $|B_{\theta+d_0+\epsilon}(z_0)\backslash E|$ before sending $\epsilon\to 0$:
\begin{equation}\frac{n-1}{\theta}+\frac{D}{\Delta t}\ge(\min\phi)^2((\theta+d_0)^{n-1}(\frac{1}{\theta^{n-2}}-\frac{1}{(\theta+d_0)^{n-2}}))^{-2}.\end{equation}
Now we use the convexity of the function $t\mapsto\frac{1}{t^{n-2}}$, and the fact $d_0<D-\theta$ to obtain $$(\theta+d_0)^{n-1}(\frac{1}{\theta^{n-2}}-\frac{1}{(\theta+d_0)^{n-2}})\le D^{n-1}(n-2)d_0/\theta^{n-3}.$$ Plug this into (3.2) we obtain $$d_0\ge C(n,\theta, D,\phi)\sqrt{\Delta t}.$$
\end{proof} 

We next give several propositions related to various density estimates. Compare with standard estimates as in other free boundary problems, our estimates hold only for small scales depending on $\Delta t$, as is typical for problems with a lower order term. 

Also, to avoid certain pathological behaviour, statements concerning radii are to be considered in the a.e. sense.

$(u, E)$ denotes a minimizer of (3.1). 
\begin{prop}
If $B_{r}(p)\subset\Dout$, then $$\frac{Per(E;B_r(p))}{r^{n-1}}\le n\omega_n+\omega_n\frac{D}{\Delta t}r.$$
\end{prop} 
\begin{proof}
Take $\tilde{E}=E\cup B_{r}(p)$, then $(u,\tilde{E})\in\mathcal{A}$.  Then optimality gives 
\begin{align*}
0&\le (Per(\tilde{E})+\frac{1}{\Delta t}\int_{\tilde{E}\Delta E_0}dist(x,\partial E_0)dx)-(Per(E)+\frac{1}{\Delta t}\int_{E\Delta E_0}dist(x,\partial E_0)dx)\\&\le (H^{n-1}(\partial B_r(p)\backslash E)-Per(E;B_r(p)))+\frac{1}{\Delta t}\int_{B_r(p)\backslash E_0}dist(x,\partial E_0)dx\\&\le (H^{n-1}(\partial B_r(p))-Per(E;B_r(p)))+\frac{D}{\Delta t}\omega_nr^n.
\end{align*}

Thus $$Per(E;B_r(p))\le n\omega_nr^{n-1}+\omega_n\frac{D}{\Delta t}r^n.$$
\end{proof} 

\begin{prop}
If $p\in\partial E$ and $B_{R_0}(p)\subset\Dout$, then 
\begin{enumerate}
\item $|E^c\cap B_{R_0}(p)|\ge c(n)R_0^n$ if $0<R_0<\frac{1}{2}n\omega_n^{\frac{1}{n}}\frac{\Delta t}{D}$, and 
\item for general $R_0$, either $|E^c\cap B_{R_0}(p)|\ge c(n)R_0^n$ or $|E^c\cap B_{R_0}(p)|\ge c(n)(\frac{\Delta t}{D})^n$.
\end{enumerate}
\end{prop} 
\begin{proof}For $0<r<R_0$, define $$m(r)=|E^c\cap B_{r}(p)|.$$
Let's take $\tilde{E}$ as in the previous lemma, then the same comparison gives

$$Per(E;B_r(p))\le H^{n-1}(\partial B_r(p)\backslash E)+\frac{1}{\Delta t}\int_{B_r(p)\backslash E}dist(x,\partial E_0)dx.$$ A usual trick in geometric measure theory gives:
$$Per(E;B_r(p))+H^{n-1}(\partial B_r(p)\backslash E)\le 2H^{n-1}(\partial B_r(p)\backslash E)+\frac{1}{\Delta t}\int_{B_r(p)\backslash E}dist(x,\partial E_0)dx.$$

Now note that the left-hand side is $Per(E^c\cap B_r(p))$, and is thus bounded from below via the isoperimetric inequality by $n\omega_n^{\frac{1}{n}}m(r)^{\frac{n-1}{n}}$. 

For the right-hand side we note $m'(r)=H^{n-1}(\partial B_r(p)\backslash E)$.

Consequently we have the following differential inequality on $m(r)$:
$$n\omega_n^{\frac{1}{n}}m(r)^{\frac{n-1}{n}}\le 2m'(r)+\frac{D}{\Delta t}m(r),$$which is equivalent to 
$$m'(r)\ge \frac{1}{2}m(r)^{\frac{n-1}{n}}(n\omega_n^{\frac{1}{n}}-\frac{D}{\Delta t}m(r)^{\frac{1}{n}}).$$

In the nice case where $\frac{D}{\Delta t}m(r)^{\frac{1}{n}}\le \frac{1}{2}n\omega_n^{\frac{1}{n}}$ for all $0<r<R_0$, one has $$m'(r)\ge \frac{1}{4}n\omega_n^{\frac{1}{n}}m(r)^{\frac{n-1}{n}}.$$ A Gronwall type estimate leads to 
$$m(R_0)\ge c(n)R_0^n.$$

In the bad case where $\frac{D}{\Delta t}m(r)^{\frac{1}{n}}\ge  \frac{1}{2}n\omega_n^{\frac{1}{n}}$ for some $r<R_0$, we simply note that $$m(R_0)\ge m(r)\ge c(n)(\frac{\Delta t}{D})^n.$$

Then note the bad case is possible only when $R_0\ge \frac{1}{2}n\omega_n^{\frac{1}{n}}\frac{\Delta t}{D}.$
\end{proof} 

Although we do not need the following directly, we show that when $x_0$ is {\it inside $E_0$}, we can improve the previous estimate with a much easier argument.

\begin{prop}
If $x_0\in\overline{E}$ and $B_r(x_0)\subset E_0$, then $$|E^c\cap B_r(x_0)|\ge c(n)r^n.$$
\end{prop} 
\begin{proof}
$(u, \tilde{E})\in\mathcal{A}$ where $\tilde{E}=E\cup B_s(x_0)$ $0<s<r$.

Since $B_s(x_0)\subset E_0$, this competitor has less distance integral and the same Dirichlet energy. Consequently optimality on $(u,E)$ implies $$Per(\tilde{E})\ge Per(E).$$

Note that $Per(\tilde{E})=H^{n-1}(\partial B_s(x_0)\cap E^c)+Per(E;B_s(x_0)^c)$, we have $$H^{n-1}(\partial B_s(x_0)\cap E^c)\ge Per(E;B_s(x_0)),$$ or equivalently, $$ 2H^{n-1}(\partial B_s(x_0)\cap E^c)\ge Per(E;B_s(x_0))+H^{n-1}(\partial B_s(x_0)\cap E^c).$$

Define $m(s)=|B_s(x_0)\backslash E|$, isoperimetric inequality and the previous estimate read $$m'(s)\ge c(n)/2m(s)^{\frac{n-1}{n}},$$ from which the desired estimate follows.
\end{proof} 

We are not particularly interested in the regularity of the function $u$ in this paper. However, the following rudimentary analysis is needed to justify certain calculations later. The argument very much follows Mazzone \cite{Mazzone} and David-Toro \cite{DT}.

\begin{prop} There is an $\alpha=\alpha(n)\in(0,1)$ such that 
for $p\in\partial E$, $B_{R_0}(p)\subset\Dout\backslash\overline{\Din}$ with $0<R_0<\frac{1}{2}n\omega_n^{\frac{1}{n}}\frac{\Delta t}{D}$, one has

$$u(x)\le \sup_{B_{R_0}(p)}u  \cdot  (\frac{1}{R_0})^{\alpha}|x-p|^{\alpha} \text{ for all $x\in B_{R_0}(x_0)$}.$$
\end{prop} 
\begin{proof}
For $0<r<R_0$ define $$v(x)=\sup_{B_r(p)}u-u(x).$$ Then $v$ is a nonnegative superharmonic function in $B_r(p)$.

By the weak Harnack inequality \cite{CC}, there is a dimensional $p_0>0$,
\begin{align*}
\inf_{B_{\frac{r}{2}}(p)}v&\ge c(n)(\frac{1}{\omega_nr^n}\int_{B_r(p)}v^{p_0})^{\frac{1}{p_0}}\\&\ge c(n)(\frac{1}{\omega_nr^n}\int_{B_r(p)\backslash E}v^{p_0})^{\frac{1}{p_0}}\\&=c(n)\sup_{B_r(p)}u(\frac{1}{|B_r(p)|}|E^c\cap B_r(p)|)^{\frac{1}{p_0}}\\&\ge c(n)\sup_{B_r(p)}u.
\end{align*}

That is, $$\sup_{B_{\frac{r}{2}}(p)}u\le (1-c(n))\sup_{B_r(p)}u.$$ 

From here the H\"older decay follows by standard iteration. 
\end{proof} 

Note that $u$ is smooth away from the free boundary. The previous lemma shows that it decays around free boundary points in a continuous fashion. Therefore,
\begin{cor}
$u$ is a continuous function. $\{u>0\}$ is open. $\Delta u=0$ in $\{u>0\}$.
\end{cor}

Next we upgrade the $\alpha$-decay to a $\frac{1}{2}$-decay.

\begin{prop}
For $x_0\in E$ with $0<d_0:=dist(x_0,\partial E)<dist(x_0,\partial\Din)$, one has $$u(x_0)\le C(n,\Dout,\Din,\phi)\sqrt{d_0+\frac{D}{\Delta t}d^2_0}.$$
\end{prop} 
\begin{proof}
Let $y_0\in\partial E$ be such that $d_0=|x_0-y_0|$.

If $y_0\in\partial\Dout$, then we compare $u$ with the solution to 
$$\begin{cases}
\Delta v=0 &\text{ in $\Dout\backslash\overline{\Din}$,}\\ v=\phi &\text{ on $\overline{\Din}$,}\\ v=0 &\text{ on $\partial\Dout$.}
\end{cases}$$ 
Elliptic regularity dictates that $v$ is smooth up to the boundary in $\Dout\backslash\overline{\Din}$. Also, since the positive phase if $u$ is contained in $\Dout$, we have $u\le v$. Thus $$u(x_0)\le v(x_0)\le C(n,\Din,\Dout,\phi)d_0^{\frac{1}{2}}.$$

Therefore we only need to consider the case $y_0\in\Dout$. In this case, we can find $\epsilon>0$ small such that $B_{d_0+\epsilon}(x_0)\subset\Dout$ and $B_{d_0+\epsilon}(x_0)\cap\overline{\Din}=\emptyset$.

Since $B_{d_0}(x_0)$ is contained in the positive phase of $u$, Corollary 3.14 tells us that $u$ is a nonnegative harmonic function in $B_{d_0}(x_0)$. As such, Harnack inequality implies $$u(y)\ge c(n)u(x_0) \text{ for all $y\in B_{\frac{d_0}{3}}(x_0)$}.$$

Define $v$ as the solution to 
$$\begin{cases}
\Delta v=0 &\text{ in $B_{d_0+\epsilon}(x_0)\backslash B_{\frac{d_0}{3}}(x_0)$},\\ v=c(n)u(x_0) &\text{ in $B_{\frac{d_0}{3}}(x_0)$},\\ v=0 &\text{ on $\partial B_{d_0+\epsilon}(x_0)$}.
\end{cases}$$

We build as a competitor $\tilde{u}=\max\{u,v\}$ and $\tilde{E}=E\cup B_{d_0+\epsilon}(x_0)$. 

With similar computations as in the proof of Proposition 3.9, we have
$$Per(\tilde{E})-Per(E)\le\frac{n-1}{d_0}|B_{d_0+\epsilon}(x_0)\backslash E|,$$and 
$$\frac{1}{\Delta t}\int_{\tilde{E}\Delta E_0}dist(x,\partial E_0)dx-\frac{1}{\Delta t}\int_{E\Delta E_0}dist(x,\partial E_0)dx\le \frac{D}{\Delta t}|B_{d_0+\epsilon}(x_0)\backslash E|.$$

Now with $v\le u$ in $B_{\frac{d_0}{3}}(x_0)$, we can also apply similar ideas as in the proof of Proposition 3.9 to the Dirichlet energy:
\begin{align*}
\int |\nabla u|^2dx-\int |\nabla\tilde{u}|^2dx &=\int_{B_{d_0+\epsilon}(x_0)\backslash B_{\frac{d_0}{3}}(x_0)}|\nabla u|^2-|\nabla\tilde{u}|^2dx\\ &\ge\int_{B_{d_0+\epsilon}(x_0)\backslash B_{\frac{d_0}{3}}(x_0)}|\nabla u|^2-|\nabla\tilde{u}|^2-2\nabla\tilde{u}\cdot\nabla(u-\tilde{u})dx\\&\ge\int_{B_{d_0+\epsilon}(x_0)\backslash E}|\nabla u|^2-|\nabla\tilde{u}|^2-2\nabla\tilde{u}\cdot\nabla(u-\tilde{u})dx\\&=\int_{B_{d_0+\epsilon}(x_0)\backslash E}|\nabla v|^2dx\\&\ge c(n)(\frac{u(x_0)}{d_0})^2|B_{d_0+\epsilon}(x_0)\backslash E|.
\end{align*}

Optimality of $(u,E)$ and $\epsilon\to 0$ gives $$(\frac{u(x_0)}{d_0})^2\le C(n)(\frac{1}{d_0}+\frac{D}{\Delta t}).$$
\end{proof} 

Harmonicity inside the positive phase and the previous decay estimate leads to the following:
\begin{cor}
For compact $K$ inside $\Dout\backslash\overline{\Din}$, we have $$\|u\|_{C^{1/2}(K)}\le C(n,\Din,\Dout,\phi,K)\sqrt{\frac{D}{\Delta t}}.$$
\end{cor}

We now give the other side of Proposition 3.11.
\begin{prop}
If $x_0\in\overline{E}$, $0<R_0<C(n)(\frac{\Delta t}{D})^{\frac{1}{n}}$ and $B_{R_0}(x_0)\subset\Dout\backslash\overline{\Din}$, then $$\frac{|E\cap B_{R_0}(x_0)|}{|B_{R_0}(x_0)|}\ge\tilde{\delta}_0=\tilde{\delta}_0(n,\phi).$$
\end{prop} 
\begin{proof}
For each $0<s<r<R_0$ define $$v(x)=\frac{|x-x_0|-s}{r-s}.$$

As a competitor define 
$$\tilde{u}=\begin{cases}
u &\text{ outside $B_r(x_0)$,}\\0 &\text{ inside $B_s(x_0)$,}\\ \min\{u,(\sup_{B_r(x_0)}u)v\} &\text{ in $B_r(x_0)\backslash B_s(x_0)$.}
\end{cases}$$
$$\tilde{E}=E\backslash B_s(x_0).$$

Then we have the following 
\begin{align*}\frac{1}{\Delta t}\int_{\tilde{E}\Delta E_0}dist(x,\partial E_0)dx-\frac{1}{\Delta t}\int_{E\Delta E_0}dist(x,\partial E_0)dx&\le \frac{1}{\Delta t}\int_{E_0\cap B_s(x_0)}dist(x,\partial E_0)dx\\&\le \omega_n\frac{D}{\Delta t}s^n.
\end{align*}

Also, $$\int|\nabla\tilde{u}|^2dx-\int|\nabla u|^2dx\le C(n)\frac{\sup_{B_r(x_0)}u^2}{(r-s)^2}|E\cap(B_r(x_0)\backslash B_s(x_0))|.$$

As a result, the optimality condition reads
\begin{equation}Per(E;B_s(x_0))-H^{n-1}(\partial B_s(x_0)\cap E)\le \omega_n\frac{D}{\Delta t}s^n+ C(n)\frac{\sup_{B_r(x_0)}u^2}{(r-s)^2}|E\cap(B_r(x_0)\backslash B_s(x_0))|.\end{equation}

Now we take $S_0=\frac{1}{2}R_0$, and for each $m\in\mathbb{N}$, $S_{m+1}=S_m-c2^{-m}R_0$, where $c$ is some constant between $1/4$ and $1/2$.

Define also $M_m=\sup_{B_{S_m}(x_0)}u^2$ and $V_m=|E\cap (B_{S_m}(x_0)\backslash B_{S_{m+1}}(x_0))|.$

By subharmonicity of $u^2$, 
$$M_{m+1}\le \sup_{B_{S_{m+1}+c2^{-m-1}R_0}(x_0)}u^2\le \sup_{\partial B_{S_{m+1}+c2^{-m-1}R_0}(x_0)}u^2.$$

Let $y_0$ be a point realizing the supremum on the right-hand side,  then $$B_{c2^{-m-1}R_0}(y_0)\subset B_{S_m}(x_0)\backslash B_{S_{m+1}}(x_0).$$

Thus again by subharmonicity of $u^2$,
\begin{align*}u^2(y_0)&\le\frac{1}{|B_{c2^{-m-1}R_0}(y_0)|}\int_{B_{c2^{-m-1}R_0}(y_0)} u^2\\&\le\frac{1}{|B_{c2^{-m-1}R_0}(y_0)|}\int_{E\cap (B_{S_m}(x_0)\backslash B_{S_{m+1}}(x_0))} u^2\\&\le\frac{1}{|B_{c2^{-m-1}R_0}(y_0)|}M_mV_m\\&=\frac{1}{\omega_n}c^{-n}2^{n(m+1)}R_0^{-n}M_mV_m.
\end{align*}That is, \begin{equation}M_{m+1}\le C(n)2^{nm}R_0^{-n}M_mV_m.\end{equation}

Now for each $r\in(S_{m+1},S_m-c2^{-m-1}R_0)$, the isoperimetric inequality implies 
\begin{align*}
V_{m+1}^{\frac{n-1}{n}}&\le |B_r(x_0)\cap E|^{\frac{n-1}{n}}\\ &\le C(n)(H^{n-1}(\partial B_r(x_0)\cap E)+Per(E;B_r(x_0)))\\ &\le C(n)(2H^{n-1}(\partial B_r(x_0)\cap E)+\omega_n\frac{D}{\Delta t}r^n+C(n)\frac{\sup_{B_{S_m}(x_0)}u^2}{(S_m-r)^2}|E\cap(B_{S_m}(x_0)\backslash B_r(x_0))|).
\end{align*}For the last inequality, we used (3.3) with $s$ replaced by $r$ and $r$ replaced by $S_m$.

With $S_m-c2^{-m-1}R_0>r>S_{m+1}$ , the expression above can be bounded
$$V_{m+1}^{\frac{n-1}{n}}\le C(n)(H^{n-1}(\partial B_r(x_0)\cap E)+\frac{D}{\Delta t}S_m^n+\frac{4^{mn}}{R_0^2}M_mV_m).$$
Integrating over the interval for $r$
\begin{align*}
V_{m+1}^{\frac{n-1}{n}}(S_m-c2^{-m-1}R_0-S_{m+1})&\le C(n)|E\cap(B_{S_m-c2^{-m-1}R_0}(x_0)\backslash B_{S_{m+1}}(x_0)|\\&+C(n)(\frac{D}{\Delta t}S_m^n+\frac{4^{mn}}{R_0^2}M_mV_m)(S_m-c2^{-m-1}R_0-S_{m+1})\\&\le C(n)(V_m+(\frac{D}{\Delta t}S_m^n+\frac{4^{mn}}{R_0^2}M_mV_m)(S_m-c2^{-m-1}R_0-S_{m+1})).
\end{align*}

Therefore \begin{align*}V_{m+1}^{\frac{n-1}{n}}&\le C(n)(\frac{2^{m+1}}{cR_0}V_m+\frac{D}{\Delta t}S_m^n+\frac{4^{mn}}{R_0^2}M_mV_m)\\&\le \frac{C(n)}{R_0^2}4^{mn}(V_m+M_mV_m+\frac{D}{\Delta t}S_m^n).\end{align*}

This inequality and (3.4) implies
$$V_{m+1}+M_{m+1}\le \frac{C(n)}{R_0^n}A^{m}(V_m+M_m+\frac{D}{\Delta t}S_m^n)^{\frac{n}{n-1}},$$where $A$ is a possibly large constant. 

By taking a larger $C(n)$ if necessary, and use the convexity of $t\mapsto t^{\frac{n}{n-1}}$, we have 
$$V_{m+1}+M_{m+1}+\frac{D}{\Delta t}S_{m+1}^n\le \frac{C(n)}{R_0^n}A^{m}(V_m+M_m+\frac{D}{\Delta t}S_m^n)^{\frac{n}{n-1}}.$$

Then it's elementary to see that there is some $\delta_0=\delta_0(n)$ such that if $V_m+M_m+\frac{D}{\Delta t}S_m^n\le\delta_0$ for some $m$, then $M_k\to 0$ as $k\to\infty$, that is $u=0$ in $B_{\frac{1}{2}R_0}(x_0)$. (Lemma 1.4.7 in Lady\u zenskaja-Ural'ceva \cite{LU}). This contradicts that $x\in\overline{E}$.

Consequently, for all $m$ $$V_m+M_m+\frac{D}{\Delta t}S_m^n\ge\delta_0.$$

With $M_m\le\sup\phi^2\frac{|E\cap B_{S_m}(x_0)|}{|B_{S_m}(x_0)|}$ by subharmonicity, we further have $$V_m+M_m\le C(n,\phi)\frac{|E\cap B_{S_m}(x_0)|}{|B_{S_m}(x_0)|}.$$

Thus $$C(n,\phi)\frac{|E\cap B_{S_m}(x_0)|}{|B_{S_m}(x_0)|}+\frac{D}{\Delta t}S^n_m\ge\delta_0.$$

In particular, if $R_0$ is small such that $\frac{D}{\Delta t}R_0^n<\frac{1}{2}\delta_0$, then 
$$\frac{|E\cap B_{S_m}(x_0)|}{|B_{S_m}(x_0)|}\ge c(n,\phi)\delta_0$$ for all $m$.  

Sending $m\to\infty$ gives the desired result. 
\end{proof}

Now we give the final piece we need before proceeding to the compactness of $\Delta t$-approximate flows. 
\begin{prop}
If $x_0\in\partial E$ $0<r\le c(n)\frac{\Delta t}{D}$ and $B_r(x_0)\subset\Dout\backslash\overline{\Din}$,  then $$c(n,\phi)\le \frac{Per(E;B_r(x_0))}{r^{n-1}}\le C(n).$$
\end{prop} 

\begin{proof}
The upper bound follows directly from Proposition 3.10.

For the lower bound, we first note that under the assumptions, we have the following estimate
$$c(n,\phi)\le \frac{|E\cap B_r(x_0)|}{r^n}\le C(n).$$ This is a consequence of Proposition 3.11 and 3.17.

Now by the isoperimetric inequality, we have the following:
$$|E\cap B_r(x_0)|^{\frac{n-1}{n}}\le\frac{1}{n\omega_n^{\frac{1}{n}}} Per(E\cap B_r(x_0)),$$ and 
$$|E^c\cap B_r(x_0)|^{\frac{n-1}{n}}\le\frac{1}{n\omega_n^{\frac{1}{n}}} Per(E^c\cap B_r(x_0)).$$
Adding these pieces 
\begin{align*}
|E\cap B_r(x_0)|^{\frac{n-1}{n}}+|E^c\cap B_r(x_0)|^{\frac{n-1}{n}}&\le\frac{1}{n\omega_n^{\frac{1}{n}}} (Per(E\cap B_r(x_0))+Per(E^c\cap B_r(x_0)))\\&=\frac{1}{n\omega_n^{\frac{1}{n}}}(Per(E;B_r(x_0))+H^{n-1}(\partial B_r(x_0)\cap E)\\&+Per(E;B_r(x_0))+H^{n-1}(\partial B_r(x_0)\cap E^c))\\&=\frac{1}{n\omega_n^{\frac{1}{n}}}H^{n-1}(\partial B_r(x_0))+\frac{2}{n\omega_n^{\frac{1}{n}}}Per(E;B_r(x_0))\\&=|B_r(x_0)|^{\frac{n-1}{n}}+\frac{2}{n\omega_n^{\frac{1}{n}}}Per(E;B_r(x_0)).
\end{align*}

Since both $|E\cap B_r(x_0)|$ and $|E^c\cap B_r(x_0)|$ are uniformly bounded away from $0$ and $|B_r|$ by the density estimates, we have the following by the concavity of $t\mapsto t^{\frac{n-1}{n}}$:
$$|E\cap B_r(x_0)|^{\frac{n-1}{n}}+|E^c\cap B_r(x_0)|^{\frac{n-1}{n}}\ge (1+c(n,\phi))|B_r|^{\frac{n-1}{n}}.$$

Therefore $\frac{2}{n\omega_n^{\frac{1}{n}}}Per(E;B_r(x_0))\ge c(n,\phi)|B_r|^{\frac{n-1}{n}}$, which gives the desired estimate. 
\end{proof}

\subsection{Existence and continuity of a flat flow}
In this subsection we give the proof of Theorem 1.2, which states the existence and H\"older continuity of a flat flow starting from any set of finite perimeter.  Here we are considering a class of very general initial configurations, the price to pay is that we lose uniqueness, the semigroup property and consistency with the smooth MCND flow. However, it is very likely that one can recover these properties for special geometries. 

For a given $E_0$, for fixed $\Delta t>0$, we have a minimizer $(u,E)$ as in the previous subsection. Denote this minimizer by $\EDt(\Delta t)$, then we can repeat the same procedure with $\EDt(\Delta t)$ as the initial configuration and obtain a minimizer $\EDt(2\Delta t)$. Iteratively we get a sequence $\{\EDt(k\Delta t)\}_{k\in\mathbb{N}}$. This is a $\Delta t$-approximate flow in the sense of Definition 3.2. The goal is compactness as $\Delta t\to 0$.

Since we never used any property of $E_0$ in the previous subsection, results in the previous subsection  hold for any set in this sequence. These estimates do depend on $\Delta t$ and blow up when $\Delta t\to 0$. As a result, it might seem that no compactness is possible. However, there is one part of the energy that behaves very well with small $\Delta t$, namely, the distance integral. 

Since total energy is non-increasing along the approximate flow, $$\frac{1}{\Delta t}\int_{\EDt((k+1)\Delta t)\Delta\EDt(k\Delta t)}dist(x,\partial E_0)dx\le C(E_0)$$ uniformly over $\Delta t$. Thus when $\Delta t$ is small, this is pushing two consecutive sets in an approximate flow to be very close to each other in measure. The brilliant idea in Almgren-Taylor-Wang is to exploit this regularization fact to cancel all the bad behaviour from other estimates. 

We would very much follow their idea. So far we have been treating the distance integral as a lower-order error, now we use crucially the uniform boundedness of this term. This begins with the following technical  lemma, which is based on Proposition 4.3 in \cite{ATW} but modified for our purpose:

\begin{lem}
Suppose that $C$ and $A$ are measurable. 

Let $\delta$, $\gamma$, $\Delta t$ and $E$ be positive numbers such that $$\frac{1}{\Delta t}\int_{A\backslash C}dist(x,\partial C)dx\le E,$$ and $$H^{n-1}(\partial C\cap B_r(p))\ge\gamma r^{n-1}$$ whenever $p\in \partial C$ and $0<r\le\delta$.  

Then for $\delta\le R<\infty$ we have $$|A\backslash C|\le[2\Gamma(\frac{R}{\delta})^{n-1}H^{n-1}(\partial C)]^{1/2}(\Delta t)^{1/2}E^{1/2}+\frac{\Delta t}{R}E,$$ where $\Gamma=2^{2n+1}n\omega_n\beta(n)/\gamma$, and $\beta(n)$ is the dimensional constant in Besicovitch covering lemma.
\end{lem} 

For the convenience of the reader, we give some ideas behind its proof. The interested reader should see \cite{ATW} for a rigorous proof.

\begin{proof}

First, let $A$ and $C$ be any set, we first prove the following lemma about distance integrals:
\begin{lem}
$$\int_{A\backslash C} dist(x,\partial C)dx\le E \implies |A\backslash C|\le 2^{1/2}[\sup_{0<r<R}H^{n-1}(A\cap\{dist(\cdot,\partial C)=r\})]^{1/2}E^{1/2}+\frac{E}{R}.$$
\end{lem}

First, by Chebychev $$|(A\backslash C)\cap \{dist(\cdot, \partial C)\ge R\}|\le \frac{E}{R}.$$ To deal with the set $|A\cap \{dist(\cdot, \partial C)\le R\}|$, we decompose it into level sets of the distance function and apply the coarea formula \cite{Maggi}:
\begin{align*}
E&\ge \int_{A\cap \{dist(\cdot, \partial C)\le R\}} dist(x,\partial C)dx\\&=\int_0^RrH^{n-1}(A\cap \{dist(\cdot, \partial C)=r\})dr.
\end{align*}
Note that $$|A\cap \{dist(\cdot, \partial C)\le R\}|=\int_0^RH^{n-1}(A\cap \{dist(\cdot, \partial C)=r\})dr,$$ it is natural to use a rearrangement argument to show 
\begin{align*}\int_0^RH^{n-1}(A\cap &\{dist(\cdot, \partial C)=r\})dr\\&\le 2^{1/2}(\int_0^RrH^{n-1}(A\cap \{dist(\cdot, \partial C)=r\})dr)^{1/2}\\&\cdot (\sup_{0<r<R}H^{n-1}(A\cap \{dist(\cdot,\partial C)=r\}))^{1/2},\end{align*} which concludes the proof for Lemma 3.20.
 
Compare this general lemma with the situation in Lemma 3.18, we simply need to deduce \begin{equation}\sup_{0<r<R}H^{n-1}(\{dist(\cdot, \partial C)=r\})\le\Gamma (\frac{R}{\delta})^{n-1}H^{n-1}(\partial C)\end{equation} from the density lower bound. 

We first establish the following fact about distance functions:
\begin{lem}
For any closed set $C$, $$H^{n-1}(\{dist(\cdot, \partial C)=1\}\cap B_2)\le 2^{2n+1}n\omega_n.$$
\end{lem} 

By cutting off, we might assume $C\subset B_3$. Then coarea formula and the mean value theorem gives some $1/2<R<1$ such that $$H^{n-1}(\{dist(\cdot, \partial C)=R \})\le 2|\{1/2\le dist(\cdot, \partial C)\le 1\}|.$$

Then we invoke \begin{align*}H^{n-1}(\{dist(\cdot, \partial C)=1\})-H^{n-1}(\{dist(\cdot, \partial C)=R\})&\le \frac{n-1}{R}|\{R\le dist(\cdot, \partial C)\le 1\}|\\&\le 3^nn\omega_n,\end{align*}which can be easily proved using the divergence theorem. 

Lemma 3.21 follows by adding the previous estimates. 

Now we turn to (3.5).

We first deal with the case when $r<\delta$. 

Cover $\partial C$ with $\{B_r(p)\}_{p\in\partial C}$, and reduce it to $\{B_r(p_j)\}$ by Besicovitch. Since $\{dist(\cdot, \partial C) =r\}$ is covered by $\{B_{2r}(p_j)\},$ \begin{align*}H^{n-1}(\{dist(\cdot, \partial C)=r\})&\le\Sigma H^{n-1}(B_{2r}(p_j)\cap \{dist(\cdot, \partial C)=r\})\\&\le\Sigma 2^{2n+1}n\omega_nr^{n-1}\\&\le 2^{2n+1}n\omega_n\Sigma \frac{1}{\gamma}H^{n-1}(\partial C\cap B_r(p_j))\\&\le \frac{\Gamma}{\gamma} H^{n-1}(\partial C).\end{align*}Note that we used a scaled version of Lemma 3.21 and the density lower bound on $\partial C$.

For the case when $r>\delta$, we simply note that $$r^{n-1}\le (\frac{r}{\delta})^{n-1}\frac{1}{\gamma}H^{n-1}(\partial C\cap B_\delta(p)) \text{ for $p\in\partial C$}.$$ After this we  can use the same covering argument but with balls of the form $B_\delta(p)$ for $p\in\partial C$.

This concludes the proof for Lemma 3.19.
\end{proof} 

The next proposition gives a uniform H\"older estimate in time for approximate flows.
\begin{prop}
Let $\{\EDt(k\Delta t)\}_{k\in\mathbb{N}}$ be a $\Delta t$-approximate flow starting from $E_0$, then for any $N\in\mathbb{N}$ one has $$|\EDt((k+N)\Delta t)\Delta \EDt(k\Delta t)|\le C(n,E_0,\phi,\Din,\Dout)(N\Delta t)^{\frac{1}{n+1}}.$$
\end{prop} 

\begin{proof}
Define $$M_k=(\int |\nabla u_{k-1}|^2dx+Per(E((k-1)\Delta t))-(\int |\nabla u_{k}|^2dx+Per(E(k\Delta t))),$$where $u_j$ is a capacity potential corresponding to $\EDt(j\Delta t)$.

By taking $(u_{k-1},\EDt((k-1)\Delta t))$ as a competitor in the energy for $(u_k,\EDt(k\Delta t))$, we see $$\frac{1}{\Delta t}\int_{\EDt(k\Delta t)\Delta\EDt((k-1)\Delta t)}dist(x,\partial \EDt((k-1)\Delta t))dx\le M_k.$$

In particular this implies 

 \begin{equation}\frac{1}{\Delta t}\int_{\EDt(k\Delta t)\backslash\EDt((k-1)\Delta t)}dist(x,\partial \EDt((k-1)\Delta t))dx\le M_k,\end{equation}and 
 $$\frac{1}{\Delta t}\int_{\EDt((k-1)\Delta t)\backslash\EDt(k\Delta t)}dist(x,\partial\EDt((k-1)\Delta t))dx\le M_k,$$ which is equivalent to 
 \begin{equation}\frac{1}{\Delta t}\int_{\EDt(k\Delta t)^c\backslash\EDt((k-1)\Delta t)^c}dist(x,\partial(\EDt((k-1)\Delta t)^c))dx\le M_k.\end{equation}

Due to Proposition 3.18, for $0<r<c(n)\frac{\Delta t}{D}$ and $x_0\in\partial\EDt((k-1)\Delta t),$ we have the density estimate
$$H^{n-1}(\partial \EDt((k-1)\Delta t)\cap B_r(x_0))\ge c(n,\phi)r^{n-1},$$ hence we can apply Lemma 3.19 to (3.6) with $A=\EDt(k\Delta t)$, $C=\EDt((k-1)\Delta t)$, $\delta=c(n)\frac{\Delta t}{D}$, $\gamma=c(n,\phi)$ and $E=M_k$ to obtain 
\begin{align*}|\EDt(k\Delta t)&\backslash\EDt((k-1)\Delta t)|\\&\le C(n,\phi)D^{\frac{n-1}{2}}(H^{n-1}(\partial \EDt((k-1)\Delta t)))^{1/2}(\frac{R}{\Delta t})^{\frac{n-1}{2}}(\Delta t)^{1/2}M_k^{1/2}+\frac{\Delta t}{R}M_k\end{align*} for all $R>c(n)\frac{\Delta t}{D}$.

Note that $Per(\EDt((k-1)\Delta t))\le C(E_0)$ by the monotonicity of energy, we have 
$$|\EDt(k\Delta t)\backslash\EDt((k-1)\Delta t)|\le C(n,\phi,D,E_0)(\frac{R}{\Delta t})^{\frac{n-1}{2}}(\Delta t)^{1/2}M_k^{1/2}+\frac{\Delta t}{R}M_k.$$

Similar arguments applied to (3.7) gives 
$$|\EDt((k-1)\Delta t)\backslash\EDt(k\Delta t)|\le C(n,\phi,D,E_0)(\frac{R}{\Delta t})^{\frac{n-1}{2}}(\Delta t)^{1/2}M_k^{1/2}+\frac{\Delta t}{R}M_k.$$
Consequently, $$|\EDt(k\Delta t)\Delta\EDt((k-1)\Delta t)|\le C(n,\phi,D,E_0)(\frac{R}{\Delta t})^{\frac{n-1}{2}}(\Delta t)^{1/2}M_k^{1/2}+2\frac{\Delta t}{R}M_k$$ for all $R>c(n)\frac{\Delta t}{D}$.

Pick $R=\frac{\Delta t}{(N\Delta t)^{\frac{1}{n+1}}}$, the previous estimate translates to 
$$|\EDt(k\Delta t)\Delta\EDt((k-1)\Delta t)|\le C(n,\phi,D,E_0)M_k^{1/2}(\frac{1}{N})^{\frac{1}{2}}(N\Delta t)^{\frac{1}{n+1}}+2M_k(N\Delta t)^{\frac{1}{n+1}}.$$

Adding up $N$ such estimates to obtain
\begin{align*}|\EDt((k+N)\Delta t)\Delta\EDt(k\Delta t)|&\le C(n,\phi,D,E_0)(N\Delta t)^{\frac{1}{n+1}}(\Sigma (M_k/N)^{\frac{1}{2}}+\Sigma M_k)\\&\le C(n,\phi,D,E_0)(N\Delta t)^{\frac{1}{n+1}}\Sigma M_k.\end{align*}

Now simply note that being a telescoping sum, $$\Sigma_{k}^{N}M_k=(\int |\nabla u_{k-1}|^2dx+Per(E((k-1)\Delta t))-(\int |\nabla u_{N}|^2dx+Per(E(N\Delta t))),$$ where the right-hand side is bounded by the energy of the initial configuration, say, $C(E_0)$.

Consequently, we have $$|\EDt((k+N)\Delta t)\Delta\EDt(k\Delta t)|\le C(n,\phi,D,E_0)(N\Delta t)^{\frac{1}{n+1}}.$$
\end{proof} 

The  proof of Theorem 1.2 follows.

\begin{proof}(of Theorem 1.2. )
For a sequence $\Delta t\to 0$, and each fixed $t>0$, one has a subsequence of $$\EDt(t)\to E(t) \text{ in $\mathcal{L}^1$}$$ for some $E(t)$. This is a consequence of the compactness of sets of finite perimeter \cite{Giusti}\cite{Maggi}. 

Then by Cantor's diagonal argument we find a subsequence that converges at all rational $t>0$. The previous uniform H\"older estimate shows the convergence happens at all real $t>0$. 

Moreover, the limiting $E:[0,\infty)\to\PR$ satisfies $|E(s)\Delta E(t)|\le C|t-s|^{\frac{1}{n+1}}.$
\end{proof} 

Note that we only have a flow at the level of the sets, not at the level of the potentials.  As a result, many interesting questions are left open about flat flows. For instance, do they satisfy the equation of motion in some weak sense? Is the energy decreasing along the flows? To tackle these questions, we need better estimate on the potentials, which might come from stronger convergence of the sets. 

Some other problems that are very challenging under this formulation of flat flows such as the uniqueness, the semigroup property and the consistency with smooth flows are easily tackled if we use the following formulation.

\section{The minimal barrier}
The continuation of the mean curvature flow by the method of minimal barriers goes back to De Giorgi \cite{DeGiorgi}. This formulation inherits a geometric comparison property of the mean curvature flow, namely, the inclusion principle. Unlike the variational approach, which is based on energy considerations, it is very easy to get comparisons in this formulation. As a result, pointwise properties like uniqueness, the semigroup property and the consistency property become almost trivial. The reader might consult Bellettini \cite{Bellettini} for more details of this approach. 

In this section, a similar formulation is given for our flow. This is based on the observation that our flow also enjoys an inclusion principle as explained in the introduction. We will see that those pointwise properties which are difficult to obtain in the variational approach will again be relatively simple under this formulation.  

But we first need some definitions. 

As the space of `test functions', we take all possible smooth MCND flows:

\begin{defi}
The class $\mathcal{F}$ consist of all smooth MCND flows in the sense of Definition 2.3.
\end{defi} 

\begin{rem}
Note that in this section we can again drop the boundedness assumption on $\Dout$ like in Section 2, since it plays no role in the analysis.
\end{rem} 
\begin{rem}
Our test functions are fixed by the boundary data $\phi$. This lack of translation invariance and reflection is one of the major difference between our theory and the theory for the mean curvature flow.
\end{rem} 

Following De Giorgi, a `supersolution' is called a barrier.
\begin{defi}
A map $\psi:[c,d]\to\PR$ is a barrier if 
\begin{enumerate}
\item $\Din\subset\subset\psi(t)\subset\subset\Dout$ for all $t$,
\item if $f:[a,b]\to\PR$ is a smooth MCND flow with $[a,b]\subset [c,d]$ and $f(a)\subset\psi(a)$, then $f(b)\subset\psi(b)$.
\end{enumerate}
The collection of barriers on $[c,d]$ is denoted by $\mathcal{B}([c,d]).$
\end{defi}

The following propositions are direct consequences of the definitions.

\begin{prop}\begin{enumerate}
\item Given any smooth closed set $\Din\subset\subset\Omega_0\subset\subset\Dout$, there is some small $T>0$ and a unique $f\in\mathcal{F}$ on $[0,T]$ with $f(0)=\Omega_0$.

\item If $f:[a,b]\to\PR\in\mathcal{F}$, then $f_t:[a+t, b+t]\to\PR$ defined by $$f_t(s)=f(s-t)$$ is also in $\mathcal{F}$. 

\item If $f:[a,b]\to\PR\in\mathcal{F}$ and $a<c<b$, then $f|_{[a,c]}$ and $f|_{[c,b]}$, the restrictions of $f$ to $[a,c]$ and $[c,b]$ respectively, are in $\mathcal{F}$.

\item If $\{\psi_{\alpha}\}_{\alpha\in\mathcal{A}}\subset\mathcal{B}[c,d]$, then $$\bigcap_{\mathcal{A}}\psi_{\alpha}\in\mathcal{B}([c,d]).$$
\end{enumerate}
\end{prop} 

Very similar to the `least supersolution' in elliptic theory, our `solution' is defined as a minimal barrier.
\begin{defi}
For $E_0\subset\R$, the minimal barrier starting from $E_0$ on $[c,d]$ is defined as 
$$\mathcal{M}(E_0)(t):=\bigcap\{\psi\in\mathcal{B}([c,d]):\psi(c)\supset E_0\}.$$
\end{defi}

The justification for this formulation is the following geometric comparison principle for smooth MCND flows.

\begin{thm}
Suppose $f,g:[a,b]\to\PR$ are two smooth MCND flowsin the sense of Definition 2.3. 

If $f(a)\subset g(a)$ then $f(b)\subset g(b)$.
\end{thm} 
\begin{proof}
To simplify notations, we denote by $d_f(\cdot, t)$ the signed distance function to $f(t)$, and by $d_g(\cdot, t)$ the signed distance function to $g(t)$. 

Let $u$ and $v$ be the potentials corresponding to $f$ and $g$ respectively.

Then it is simple to see $$f(t)\subset g(t) \text{ if and only if } d_f(\cdot, t)\ge d_g(\cdot, t).$$

Since $f(a)\subset g(a)$, the following is well-defined $$\overline{t}=\sup_{f(t)\subset g(t)} t.$$

Let $t_n$ be a sequence that converge to $\overline{t}$ with $f(t_n)\subset g(t_n)$ for each $n$.  Then for any $x$, one has by continuity $$d_f(x,\overline{t})-d_g(x,\overline{t})=\lim(d_f(x,t_n)-d_g(x,t_n))\ge 0,$$which is another way of saying $$f(\overline{t})\subset g(\overline{t}).$$ 

We next show that if $\overline{t}<b$, then we can find some $\delta>0$ such that $$f(\overline{t}+\delta)\subset g(\overline{t}+\delta),$$ leading to a contradiction to the definition of $\overline{t}$.

First, at points $d_f(x,\overline{t})>d_g(x,\overline{t})$, by continuity we can find some $\delta_x>0$ such that the inequality remains true over $[\overline{t},\overline{t}+\delta_x]$.

Now at some point where $d_f(x,\overline{t})=d_g(x,\overline{t})\ge 0$, we first find $y\in\partial f(\overline{t})$ with $d_f(x,\overline{t})=|x-y|$.

With $f(\overline{t})\subset g(\overline{t})$, we see $y\in g(\overline{t})$. If $y\in Int(g(\overline{t}))$, then 
$$d_g(x,\overline{t})<|x-y|=d_f(x,\overline{t}),$$ leading to a contradiction. Thus $y\in\partial g(\overline{t})$.

Thus we have $d_f(\cdot,\overline{t})\ge d_g(\cdot,\overline{t})$ and $d_f(y,\overline{t})=0=d_g(y,\overline{t})$, therefore $\Delta_x d_f\ge\Delta_x d_g$ at $(y,\overline{t})$.

That is, $H_{\partial f(\overline{t})}(y)\ge H_{\partial g(\overline{t})}(y)$.

Meanwhile, $f(\overline{t})\subset g(\overline{t})$ induces $0\le u< v$. With $u(y)=0=v(y)$, one has $$0\ge\unv(y,\overline{t})> v_{\nu}(y,\overline{t}).$$ 

Combining these estimates, one has $$V_{\partial f(\overline{t})}(y)< V_{\partial g(\overline{t})}(y),$$where $V$ is the outward normal velocity. 

Then by the definition of a MCND flow, $$\frac{\partial}{\partial t}d_f(x,\overline{t})>\frac{\partial}{\partial t}d_g(x,\overline{t}),$$ thus again we find some $\delta_x$ such that $d_f(x,t)\ge d_g(x,t)$ for $t\in [\overline{t},\overline{t}+\delta_x]$.

Similar argument applies to the case when $d_f(x,\overline{t})=d_g(x,\overline{t})< 0.$

Consequently, we find for every point $x$ some $\delta_x>0$ so that the inequality between $d_f$ and $d_g$ remains true on $[\overline{t},\overline{t}+\delta_x]$. By compactness of $g(\overline{t})$, this $\delta_x$ can be chosen uniformly for all $x\in g(\overline{t})$, which is enough to conclude $$f(\overline{t}+\delta)\subset g(\overline{t}+\delta)$$ for some $\delta>0$. 
\end{proof} 
\begin{rem}
This is saying $\mathcal{F}\subset\mathcal{B}$ on the same time interval.
\end{rem} 
We now prove some properties of the flow of minimal barriers. These properties actually hold for very general minimal barrier flows. See Bellettini \cite{Bellettini}.

The first proposition states that the minimal barrier starting from $E_0$ takes $E_0$ as its initial value.

\begin{prop}For the minimal barrier starting from $E_0$ on $[0,\delta]$, 
$\mathcal{M}(E_0)(0)=E_0$.
\end{prop} 

\begin{proof}
It follows from definition that $$E_0\subset \mathcal{M}(E_0)(0).$$

Now define a map $\psi:[0,\delta]\to\PR$ by $0\mapsto E_0$ and $t\mapsto \mathcal{M}(E_0)(t)$. We show that $\psi\in\mathcal{B}([0,\delta])$.

To see this, take $f\in\mathcal{F}$ on $[a,b]\subset[0,\delta]$, and $f(a)\subset\psi(a)$.

If $a=0$, then $f(0)\subset E_0$, and $f(0)$ is contained at time zero in any barrier used in the definition of the minimal barrier. Consequently, $f(b)$ is contained in any barrier in the definition of the minimal barrier at time $b$. This implies $f(b)\subset\mathcal{M}(E_0)(b)=\psi(b)$.

If $a>0$, then from item 4 in Proposition 4.5 one has $f(b)\subset\psi(b)$. We conclude from here that $\psi\in\mathcal{B}([0,\delta])$.

As a result, $\psi$ is admissible in the definition of a minimal barrier. Thus $\psi\supset\mathcal{M}(E_0)$. This inclusion at $t=0$ implies $$E_0\supset\mathcal{M}(E_0)(0).$$

\end{proof} 

In contrast to the flat flow, we have a uniqueness property, which follows directly from the definition:

\begin{prop}
Given $E_0$, the flow of minimal barriers is unique. 
\end{prop} 

It also follows from definition that this flow respects set inclusion:
\begin{prop}
If $E_0\subset F_0$, then $\mathcal{M}(E_0)(t)\subset\mathcal{M}(F_0)(t)$.
\end{prop} 

Barriers expand faster than smooth flows. But the minimal barrier expands also slower than smooth flows:
\begin{prop}
Let $f:[a,b]\to\PR$ be a smooth MCND flow. If $E_0\subset f(a)$, then $$\mathcal{M}(E_0)(b)\subset f(b).$$
\end{prop} 
\begin{proof}
Since $f$ itself is a barrier defining the minimal barrier, the conclusion follows directly from the definition. 
\end{proof}

Also, by simple comparison it is easy to prove the semigroup property. Again this is very different from the flat flows:
\begin{thm}For $t_1<t_2$,
$$\mathcal{M}(E)(t_2)=\mathcal{M}(\mathcal{M}(E)(t_1))(t_2-t_1).$$
\end{thm} 

\begin{proof}
Define a map $\psi:[0,t_2]\to\PR$ by $$\psi(t)=\begin{cases}\mathcal{M}(E)(t) &\text{ for $t\in [0,t_1]$,}\\\mathcal{M}(\mathcal{M}(E)(t_1))(t-t_1) &\text{ for $t\in[t_1,t_2]$.}\end{cases}$$
Note that Proposition 4.7 ensures that there is no ambiguity at the time $t_1$.

We show that $\psi$ is a barrier. 

Again, from item (4) in Proposition 4.4, we see $\psi\in\mathcal{B}([0,t_1])$ and $\psi\in\mathcal{B}([t_1,t_2])$. Thus it suffices to check the barrier condition for some $f\in\mathcal{F}$ on $[a,b]$ with $a\in [0,t_1)$ and $b\in (t_1,t_2]$. To this end, we note that if $$f(a)\subset\psi(a)=\mathcal{M}(E)(a),$$ then $$f(t_1)\subset\mathcal{M}(E)(t_1).$$ But this ensures $\psi(b)\supset f(b)$ since $\mathcal{M}(\mathcal{M}(E)(t_1))$ is a barrier on $[t_1,t_2]$.

We conclude that $\psi$ is a barrier and thus $\psi\supset\mathcal{M}(E)$. At $t_2$, $$\mathcal{M}(E)(t_2)\subset\mathcal{M}(\mathcal{M}(E)(t_1))(t_2-t_1).$$

The other direction follows from the fact that $\mathcal{M}(E)$ is a barrier and it takes the value $\mathcal{M}(E)(t_1)$ at $t_1$.
\end{proof} 

The next property we prove for the flow of minimal barriers is the consistency with smooth flows. Again compare this with the case for flat flows:
\begin{thm}
Let $f:[a,b]\to\PR$ be a smooth MCND flow. Then $$\mathcal{M}(f(a))(t)=f(t)$$ for $t\in[a,b]$.
\end{thm} 

\begin{proof}
Proposition 4.12 gives $f(t)\supset \mathcal{M}(f(a))(t)$.

For the other direction, note that for any $\psi\in\mathcal{B}([a,b])$ with $\psi(a)\supset f(a)$, we have $\psi(t)\supset f(t)$. 

Since $\mathcal{M}$ is taken to be the intersection of such maps, $$\mathcal{M}(f(a))(t)\supset f(t).$$
\end{proof} 

The last property is a conditional result concerning the long-term behaviour of our flow. It says if one has that starting from parallel surfaces of $\partial\Din$, smooth flows exist globally in time and converge to the optimal configuration to the elliptic problem, then so does the flow of minimal barriers.

\begin{thm}
Let $\Dout=\R$, and define for $r>0$, $D_r:=\{x\in\R:d(x,\Din)\le r\}.$

Suppose the smooth MCND flow starting from $D_r$ exists globally in time and as $t\to\infty$ converges in the Hausorff distance to the optimal configuration in the elliptic problem. 

Then starting from any bounded $E_0$ that contains $\Din$ compactly, the flow of minimal barriers starting from $E_0$ converges to the optimal configuration.
\end{thm} 

\begin{proof}
Simply note that such a initial configuration $E_0$ is trapped for some $r_1<r_2$ $D_{r_1}\subset E_0\subset D_{r_2}$. Then our flow of minimal barriers is trapped between the smooth flows starting from these sets. 
\end{proof} 

\begin{rem}
Note that the assumptions are satisfied for the radial example in Section 2. Thus starting from any initial data, the flow of minimal barriers converges to $B_{R_{opt}}$. 

\end{rem}

\section*{Acknowledgements}
The author would like to thank his PhD advisor Prof. Luis Caffarelli for his constant encouragement and guidance. The author is also grateful to Xavier Ros-Oton, Yijing Wu and Wen Yang for many fruitful discussions.



\end{document}